\documentclass[AMA, STIX1COL,12pt]{WileyNJD-v2}

\articletype{Original research article}%

\received{}
\revised{}
\accepted{}

\usepackage{moreverb}
\usepackage{fancyhdr}                         
\usepackage{psfrag}
\usepackage{color}  
\usepackage{usrextra}                         
\usepackage{lscape}
\usepackage{amsmath,exscale}
\usepackage{bbm}  
\usepackage{caption}
\usepackage{subcaption}
\usepackage{hd}                         

\begin{document}

\title{Growth-induced instabilities for transversely isotropic hyperelastic materials}

\author[1]{Cem Altun}
\author[2]{Ercan G\"{u}rses}
\author[1]{H\"usn\"u Dal*}

\authormark{Altun \textsc{et al}}

\address[1]{\orgdiv{Department of Mechanical Engineering}, \orgname{Middle East Technical University}, \orgaddress{\state{Ankara}, \country{T\"urkiye}}}

\address[2]{\orgdiv{Department of Aerospace Engineering}, \orgname{Middle East Technical University}, \orgaddress{\state{Ankara}, \country{T\"urkiye}}}

\corres{*H\"usn\"u Dal, \email{dal@metu.edu.tr}}

\presentaddress{Dumlup\i nar  Bulvar\i~1, Makina M\"uhendisligi E-203, 06800 \c{C}ankaya, Ankara, T\"urkiye}

\abstract[Summary]{This work focuses on planar growth-induced instabilities in three-dimensional bilayer structures, i.e., thick stiff film on a compliant substrate. Growth-induced instabilities are examined for a different range of fiber stiffness with a five-field Hu-Washizu type mixed variational formulation. The quasi-incompressible and quasi-inextensible limits of transversely isotropic materials were considered. A numerical example was solved by implementing the \emph{T2P0F0} element on an automated differential equation solver platform, \emph{FEniCS}. It was shown that both the wavelength and critical growth parameter $g$ decrease by increasing the fiber stiffness for the first instability, which is obtained along the stiff fiber direction. The effect of the fiber stiffness is minor on the secondary buckling, which was observed perpendicular to the fiber direction. For a range of fiber stiffnesses, bifurcation points of instabilities were also determined by monitoring displacements and energies. The energy contributions of layers with different ranges of fiber stiffnesses were examined. It is concluded that the energy release mechanism at the initiation of the primary buckling is mainly due to isotropic and anisotropic contributions of the stiff film layer. For high fiber stiffnesses, the effect of the anisotropic energy on the first buckling becomes more dominant over other types. However, in the secondary instability, the isotropic energy of the film layer becomes the dominant one. Numerical outcomes of this study will help to understand the fiber stiffness effect on the buckling and post-buckling behavior of bilayer systems.}


\maketitle

\section{Introduction}\label{sec01}
Growth-induced deformations are common phenomena that confront living systems and engineering applications. In a certain level of growth, called critical growth, structural instabilities may rise as an indicator of pathological conditions of biological tissues, and it plays a key role in the development of materials for biomedical applications. These materials are mostly in the form of multi-layer structures, which may have different mechanical properties and fiber reinforcement layouts. The main mechanism of the buckling in a multi-layered structure is the compression type of loading that causes different levels of stresses in the thin stiff film and the compliant substrate; thus, buckling leads to a release of the energy. Fiber stiffness is a critical parameter that designates the critical growth in primary buckling as a wrinkle form. It is significant to understand the mechanism behind the growth-induced instabilities for transversely isotropic materials and to observe the effect of fiber stiffness over primary and secondary instabilities on bilayer systems. Growth-induced deformations may lead to instabilities as they turn into different patterns on advanced engineering materials and soft biological tissues. Mechanism of growth is also used to achieve to desired geometrical form change in biomedical applications. Growth phenomena are widely shown up in nature as in living organisms (plants, tissue, etc.)  and in engineering materials such as polymeric gels \cite{Mora2006} and stretchable electronics \cite{Khang2009}.  A wide range of studies have been performed to understand the underlying mechanics of the tissue development subjected to growth and growth-induced instabilities and their biomedical applications such as; folding of the airway \cite{Eskandari2013}, cortical folding of the brain \cite{Budday2014}, cardiac growth \cite{Ulerich2010,Goktepe2010,Goktepe2011,Raush2010}, wrinkling on spherical geometries \cite{Li2011,Javili2014}, wrinkling on the skin \cite{Genzer2006,Tepole2011}, artery growth \cite{Larry1995, Hariton2007, Demirkoparan2007, Volokh2008, Saez2014}, growing mechanics of muscle \cite{Lubarda2002},  instabilities on thin stiff film on a compliant substrate \cite{Huang2005,Brau2011,Cai2011,Jia2013}, morphogenesis of the plates \cite{Dervaux2009,Li2022a}, buckling of swelling hydrogels \cite{Mora2006,Ionov2013,Liu2015, Liu+Zhang2015, Nardinocchi2015, Nardinocchi20152, Liu2016, Wang2021} and torsional actuator modeling \cite{Liu2003, Fang2011}, to mention but a few.  For additional information, we refer to the state-of-the-art reviews  on growth, remodelling and morphogenesis of biological tissues \cite{Larry1995, Kuhl2003, Li2012, Ambrosi2011, Menzel2012, Kuhl2014, Goriely2018}.

	Most living systems are composed of multi-layered structures  possessing different mechanical properties as in biological tissues such as artery, airwall, skin, etc. Such structures are also encountered in flexible electronics \cite{Khang2009}. In nature, most living systems have one or two fiber family reinforcement embedded in different layers of the system for different purposes under various conditions. 

	Instabilities of isotropic bilayer plates composed of film and bonded substrate were examined experimentally \cite{Khang2009, Brau2011, Yang2010}, analytically, and numerically \cite{Mora2006,Huang2005, Cai2011, Jia2013,Audoly2008p3, Javili2015, Dortdivanoglu2017}. Khang et al. \cite{Khang2009} experimentally observed the mechanical buckling of a micro-scale bilayer system; herewith, they have ended up with thickness, wavelength, strain, and delamination relations. They also highlighted that these instabilities could be used as a measure for elastic moduli of materials. Brau et al. \cite{Brau2011} studied the energy requirements of membranes to obtain a new kind of instability initiated as period-doubling both theoretically and numerically. Huang et al. \cite{Huang2005} worked on wrinkles in a layered structure experimentally and numerically in two dimensions. They obtained a relationship between the wavelength and the amplitude of wrinkles for substrates with various elastic moduli and thickness; they observed different instability forms such as stripes, labyrinths, and herringbones. In the series works \cite{Audoly2008p3, Audoly2008p1,Audoly2008p2}, different patterns (stripes, varicose, checkerboard, hexagonal) were analyzed for a bilayer plate under biaxial residual compressive stress in the film. Then the formation of herringbone (zigzag) pattern with increasing residual stress was studied numerically based on a simplified buckling model. Moreover, in the third part of the study,  an asymptotic solution was proposed for a plate with an elastic foundation in the limit of large strains. Regarding energy contributions of different instability patterns, Cai et al. \cite{Cai2011} focused on a bilayer film/substrate structure subjected to equi-biaxial compressive stress numerically and experimentally, where they observed checkerboard, hexagonal, and herringbone patterns. They ranked the modes in terms of energy, and they ended up with the fact that herringbone mode has the lowest energy, then come the checkerboard with hexagonal and triangular modes. On the other hand, Jia \cite{Jia2013} showed that hexagonal patterns minimize the elastic energy, which is the dominant mode. Javili et al. \cite{Javili2015} studied growth-induced instabilities based on eigenvalue analysis which does not impose perturbation-dependent definitions. They tested their approach in slender beam and growing film on a soft substrate, and as its outcome, an objective solution was reached. Since this method provides a good initial guess for a nonlinear analysis, it cannot be used for post-buckling behavior. As a follow-up work, Dortdivanlioglu et al. \cite{Dortdivanoglu2017} proposed isogeometric analysis (IGA) enhanced with eigenvalue analysis for a thin stiff film on a compliant substrate subjected to compressive stress. They also compared the isogeometric analysis with finite element analysis, resulting in more accurate predictions with IGA. The relationship between film thickness, film and substrate stiffness ratio, number of wrinkles, and wavelengths were also studied within this study. Additional physical effects were considered by Alawiye et al. \cite{Alawiye2019}. In \cite{Alawiye2019}, linear analysis for wrinkling problems with additional considerations of pressure, surface tension, an upper substrate, and fibers were performed. Diffusion-driven time transient swelling on hydrogels was studied in \cite{Ilsend2019, Dortdivanoglu2019}. Considering the numerical efficiency, Kadapa et al. \cite{Kadapa2020} proposed a finite element framework by extending a mixed displacement-pressure formulation using quadratic and linear Bezier elements. They concluded that the model provides good accuracy for exact incompressibility by providing an inf-sup stability condition.

	Besides isotropic modeling of bilayer plate, there are also studies on multi-layer structural deformations and instabilities taking into account the fiber reinforcement on layers. These types of structures have mostly been revealed in living systems (bone, muscle, arteries, airway, heart etc.) and biomedical engineering designs such as actuators. As one of the first sets of studies, Rachev et al. \cite{Rachev1996} studied over remodeling caused by blood pressure of arterial wall. Taber \cite{Larry1998} also proposed a growth law for arteries, including orthotropic layers (intima/media and adventitia) reinforced by muscle fibers. Lubarda et al. \cite{Lubarda2002} presented a constitutive theory for stress-dependent evolution equations for isotropic, transversely isotropic, and orthotropic biomaterials. Then, Ciarletta et al. \cite{Ciarletta2012} worked on the combined effect of the growth and anisotropy in the epithelial formation using mixed polar coordinates. By the numerical study, they concluded that the result that distribution of residual strains and mechanical properties of fibers embedded in the tissue has a significant effect on instability patterns. Liu et al. \cite{Liu2016} proposed a nonlinear finite element procedure for anisotropic swelling that included the chemical potential for anisotropic hydrogel-based bilayers. Although they do not focus on instabilities, they obtained the effects of modulus and fiber orientations in free bending and twisting. Stewart et al. \cite{Stewart2016} worked on wrinkling instability of bilayer cerebral cortex (grey and white matter) tissue embedded with two family elastic fibers in two dimensions. They showed that wrinkle wavelength is a function of fiber orientation, and the instability can be triggered by increasing the fiber stiffness depending on the fiber angle. There are also studies for macroscopic instabilities taking into account micromechanics based on homogenization \cite{Agoras2009}, Michel et al. \cite{Michel2010} extended to microstructural and macroscopic instabilities for plane-strain problems, and then Slesarenkoa et al. \cite{Slesarenkzdunek2014017} focused on macroscopic and microscopic instabilities in three-dimensional periodic fiber-reinforced composites.  They found that  the volume fraction of fibers determines the first mode of buckling and a high fraction ratio of fibers exceeding a threshold value ends up with long wave instability, whereas a lower fraction of fibers results in microscopic instabilities. Growth-induced nonlinear behavior such as remodeling of fibers \cite{Larry1995, Hariton2007, Ambrosi2011, Menzel2012, Topol2015, Lanir2015}, free swelling of anisotropic hydrogels \cite{Liu+Zhang2015,Nardinocchi2015, Nardinocchi20152, Wang2021}, anisotropic growth of the heart \cite{Goktepe2011, Raush2010}, swelling of tracheal angioedema \cite{Gou2016, Gou2017, Gou2018} and growth in fiber-reinforced torsional actuators \cite{Liu2003, Fang2011} are extensively studied. For a detailed review accounting for growth-induced nonlinear behavior of fiber-reinfored systems, the interested reader is referred to \cite{Larry1995, Ambrosi2011, Menzel2012, Kuhl2014}. 
\textsl{Albeit tremendous amount of work devoted to the aforementioned aspects of growth and remodelling, there is still need for research regarding the  the effect of fiber stiffness on growth for primary and secondary instabilities in three-dimensional bilayer fiber-reinforced structures.}

Elastomers, hydrogels, and soft tissues  exihibit nearly incompressible hyperelastic mechanical behavior that can be characterized by a free energy density function. Fiber-reinforced soft matrix materials and biological tissues in the inextensibility limit exhibit nearly incompressible and inextensible mechanical reponse. Thus, development of efficient and robust finite element formulations in the quasi-inextensible and -incompressible limit becomes crucial. Hyperelastic materials exhibit stiff volumetric response compared to shear response caused by nearly incompressible behavior. Standard displacement-based formulations show poor convergence behavior and inaccurate results for incompressible materials. A similar problem has been revealed in the nearly inextensible limit due to the high stiffness in the fiber directions \cite{zdunek2014, zdunek2016, Dal2018}. Mixed or hybrid element formulations based on variational formulations utilize additional independent variables such as stress or strain as Lagrange multipliers. In this context, we refer to the pioneering works of Pian et al. \cite{pian1982,pian184}, which were based on the Hellinger-Reissner formulation. These formulations improve the stress approximation of the standard displacement formulation. The mean dilatation formulation with \emph{Q1P0}-element based on Hu-Washizu type variational principle was introduced by Nagtegaal et al. \cite{nagtegaal1974} and extended to large-strain problems by Simo et al. \cite{simo1985}. It was implemented to hyperelastic materials within a nearly incompressible limit by Simo and Taylor \cite{simo1991}. In the element formulation, an additional term was embedded to the potential function as a constraint for the incompressibility \cite{miehe1994}, adopted in finite element implementation for visco-elastic materials \cite{dal2009}, and a novel two-field mixed displacement-pressure formulation was presented by \cite{Kadapa2020b} that provides consistent results to three-field formulation. The use of dilatation formulation for fiber-reinforced rubberlike materials and fiber-reinforced soft tissues was presented in \cite{weiss1996}.

A mixed variational formulation based on the Hu-Washizu principle for the nearly incompressible and nearly inextensible limits for fiber-reinforced materials was studied in \cite{zdunek2014,zdunek2016,schroder2016,wriggers2016,Dal2018} for one family fibers. Later it was extended to two family fibers for soft biological tissues \cite{osman2020}. The approach of Zdunek et al. \cite{zdunek2014,zdunek2016} is based on multiplicative decomposition of the deformation gradient into a purely unimodular extensional part, a purely spherical part, and an extension free unimodular part. The mixed element formulation of  Dal et al.\cite{Dal2018, osman2020} is based on a five-field Hu-Washizu type variational formulation that has conjugate pairs $(p,\theta)$ and $(s,\lambda)$ for pressure-dilatation and fiber stress-fiber stretch, respectively. The result of the variational formulation is the \emph{Q1P0F0} element which is extended for the inextensibility limit. For a detailed review, we refer to Dal et al.\cite{Dal2018,osman2020} and references therein.

In this study, we investigate three-dimensional, transeversely isotropic bilayer plates with a stiff thin film adhered to a compliant substrate under planar growth by using open-source, Python-based automated finite element platform \emph{FEniCS}  \cite{Logg2012}. It is also aimed to observe the effects of fiber stiffness over instabilities with the computational framework of the extended five-field Hu-Washizu type variational formulation proposed in \cite{Dal2018} for growth problems using the \emph{T2P0F0} element formulation. It is essential to capture and understand the mechanics behind critical instability conditions in three dimensions and to create relations between the fiber-reinforced bilayer system and its behavior under growth conditions. The paper is organized as follows: In \emph{Section} \ref{sec02}, field equations and corresponding state variables for fiber-reinforced transversely isotropic hyperelastic solid will be presented for growth-induced problems. Corresponding stress and moduli expressions are derived for both Lagrangian and Eulerian configurations. In \emph{Section} \ref{sec03}, the mixed variational formulation, which leads to the quasi-incompressible and quasi-inextensible element formulation, will be introduced together with the definition of growth tensor. In \emph{Section} \ref{sec04}, growth-induced instabilities in a three-dimensional bilayer system with one family fiber reinforcement are numerically solved with a five-field Hu-Washizu type \emph{T2P0F0} element formulation under planar growth. Fiber stiffness effect over the buckling and critical growth parameters are investigated. The results are summarized in \emph{Section} \ref{sec05}.
\section{Governing equations of motion} \label{sec02}
This section presents the kinematics and governing field equations of growth phenomena in the transversely isotropic hyperelastic solids. The kinematics based on the multiplicative decomposition of the deformation gradient into elastic and growth tensors in the sense of Rodriguez \cite{Rodriguez1994} are introduced. A simple neo-Hookean hyperelastic response in the intermediate state is considered. This study presents a quasi-incompressible and quasi-inextensible transversely isotropic model with transversely isotropic or planar growth. Finally, a five-field mixed variational formulation in the sense of Dal \cite{Dal2018} and respective Euler-Lagrange  equations are demonstrated. The governing equations are presented both in Lagrangian and Eulerian configurations. It is concluded that pursuing with the Eulerian setting leads to a more straightforward and compact form in terms of stress and moduli expressions.
\subsection{Geometric mappings and the field variables}\label{sec02-1}
\begin{figure}[h]
\centering
\def\svgwidth{0.8\textwidth}
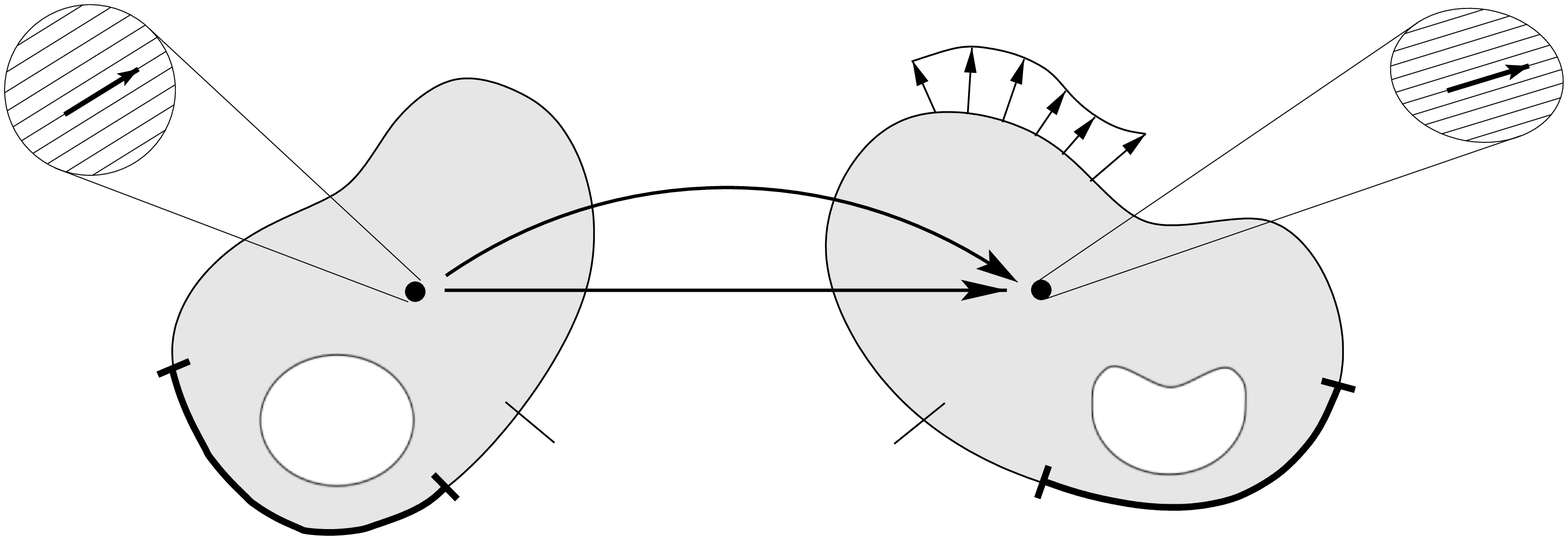
\caption{Nonlinear deformation of a transversely isotropic solid. The reference configuration $\SB\in\IR^3$ and the spatial configuration $\SSa\in\IR^3$. $\Bvarphi: \SB\times\IR \mapsto \IR^3$ is the nonlinear placement field which maps at time $t \in \IR_+$ material point position $\BX \in \SB$ onto spatial position $\Bx = \Bvarphi(\BX, t) \in \SSa$. The deformation gradient $\BF$ maps a Lagrangian infinitesimal line element $\textrm{d}\BX $ onto its Eulerian counterpart $\textrm{d}\Bx = \BF \textrm{d}\BX $.}
\label{deformation-map}
\end{figure}
\begin{figure}[h!]
\centering
\def\svgwidth{0.4\textwidth}
\begingroup%
  \makeatletter%
  \providecommand\color[2][]{%
    \errmessage{(Inkscape) Color is used for the text in Inkscape, but the package 'color.sty' is not loaded}%
    \renewcommand\color[2][]{}%
  }%
  \providecommand\transparent[1]{%
    \errmessage{(Inkscape) Transparency is used (non-zero) for the text in Inkscape, but the package 'transparent.sty' is not loaded}%
    \renewcommand\transparent[1]{}%
  }%
  \providecommand\rotatebox[2]{#2}%
  \newcommand*\fsize{\dimexpr\f@size pt\relax}%
  \newcommand*\lineheight[1]{\fontsize{\fsize}{#1\fsize}\selectfont}%
  \ifx\svgwidth\undefined%
    \setlength{\unitlength}{553.438961bp}%
    \ifx\svgscale\undefined%
      \relax%
    \else%
      \setlength{\unitlength}{\unitlength * \real{\svgscale}}%
    \fi%
  \else%
    \setlength{\unitlength}{\svgwidth}%
  \fi%
  \global\let\svgwidth\undefined%
  \global\let\svgscale\undefined%
  \makeatother%
  \begin{picture}(1,0.67249013)%
    \lineheight{1}%
    \setlength\tabcolsep{0pt}%
    \put(0,0){\includegraphics[width=\unitlength]{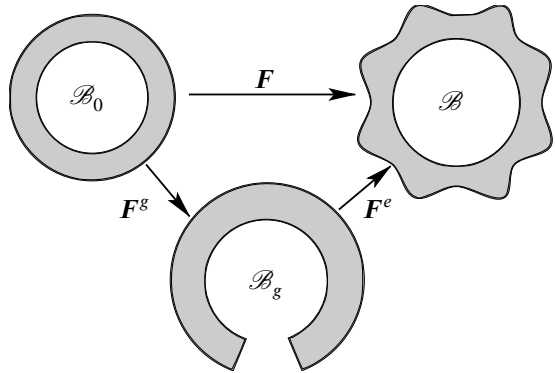}}%
    \put(0.45,0.52){\color[rgb]{0,0,0}\makebox(0,0)[lt]{\lineheight{1.25}\smash{\begin{tabular}[t]{l}$\BF$\end{tabular}}}}%
    \put(0.20,0.28){\color[rgb]{0,0,0}\makebox(0,0)[lt]{\lineheight{1.25}\smash{\begin{tabular}[t]{l}$\BF^g$\end{tabular}}}}%
    \put(0.65,0.28){\color[rgb]{0,0,0}\makebox(0,0)[lt]{\lineheight{1.25}\smash{\begin{tabular}[t]{l}$\BF^e$\end{tabular}}}}%
    \put(0.11,0.49){\color[rgb]{0,0,0}\makebox(0,0)[lt]{\lineheight{1.25}\smash{\begin{tabular}[t]{l}$\SB_0$\end{tabular}}}}%
    \put(0.79,0.48){\color[rgb]{0,0,0}\makebox(0,0)[lt]{\lineheight{1.25}\smash{\begin{tabular}[t]{l}$\SB$\end{tabular}}}}%
    \put(0.44,0.15){\color[rgb]{0,0,0}\makebox(0,0)[lt]{\lineheight{1.25}\smash{\begin{tabular}[t]{l}$\SB_g$\end{tabular}}}}%
  \end{picture}%
\endgroup%

\caption{Kinematic representation of finite growth. The multiplicative decomposition of the deformation gradient  $\BF$ gives mapping relation based on growth state $\BF^g$, and the elastic state of the deformation gradient $\BF^e$. There are defined three configurations at finite growth, the first one is the original stress-free configuration in $\SB_0$, the second one is the growth state with stress-free intermediate configuration in $\SB_g$ which also leads to incompatibility in general, and the third is the stressed-state in deformed configuration,$\SB$.}
\label{multiplicativedeformation-map}
\end{figure}

 A solid \emph{body} $\SB$ is a three-dimensional manifold, which including material points $\SP\in \SB$. 
The motion of the body can be described by a one-parameter function of time via bijective mappings
\begin{equation}\label{eq1a}
\Bchi(\SP,t)= \left\{ 
\begin{array}{rcl}
\SB   &\rightarrow & \SB(\SP,t) \in \IR^3 \times \IR_+ \\
\SP &\mapsto & \Bx=\Bchi_t(\SP)=\Bchi(\SP,t). 
\end{array}\right.
\end{equation}
The point $\Bx=\Bchi(\SP,t)$ implies the spatial configuration of the particle $\SP$ at time $t\in \IR_+$.  Let the reference configuration of the material points at a reference time $t_0$ be expressed by  $\BX=\Bchi(\SP,t_0)\in\IR^3$ and the configuration at time frame $t$ denoted by $\Bchi_t(\SP)=\Bchi(\SP,t)$. The placement map $\Bvarphi_t=\Bchi_t\circ\Bchi_{t_0}^{-1}(\BX)$ such that
\begin{equation}
\Bvarphi_t(\BX)= \left\{ 
\begin{array}{rcl}
\SB_0   &\rightarrow & \SB \in \IR^3  \\
\BX &\mapsto & \Bx=\Bvarphi(\BX,t)
\end{array}\right.
\end{equation}
maps the reference configuration $\BX\in\SB_0$ of a material point onto the spatial counterpart $\Bx\in\SB$, see Figure \ref{deformation-map}.
The \emph{deformation gradient}
\begin{equation}
 \BF:T_X\SB_0\rightarrow T_x\SB ;\quad \BF := \nablaX{\phit(\BX) =\BF^e \BF^g}
\end{equation}
maps the unit tangent of the reference or the \emph{Lagrangian} configuration onto its counterpart in the spatial or the \emph{Eulerian} configuration.  The gradient operators $\nablaX[\bullet]$ and $\nablax[\bullet]$ express the spatial derivative with respect to the reference $\BX$ and current $\Bx$ coordinates, respectively. In the finite growth formulation, 
the key kinematic definition is the multiplicative decomposition of the deformation gradient $\BF$ into a reversible elastic part $\BF^e$, and an irreversible growth part $\BF^g$ \cite{Rodriguez1994}, see Figure \ref{multiplicativedeformation-map}.
Let $d\BX$, $d\BA$, and $ dV$ represent the infinitesimal line, area, and volume elements, in the undeformed configuration. Then, the deformation gradient $\BF$, its cofactor $\cof [\BF]= \det [\BF]\BF^{-T}$ and the Jacobian $J:=\det [\BF]>0$ characterize the deformation of infinitesimal line, area, and volume elements
\eb
d\Bx = \BF d\BX~, \qquad d\Ba = \cof [\BF] d\BA~, \qquad dv= J dV~.
\ee
In the same way as split of the deformation gradient $\BF$, the Jacobian $J$ is also multiplicatively decomposed into reversible elastic volume change $J^e=\det [\BF^e]$ and an irreversible grown volume change $J^g=\det [\BF^g]$\cite{Raush2014}.
\eb
J:=\det [\BF]=J^eJ^g
\ee

\begin{figure}
\begin{center}
\def\svgwidth{0.5\textwidth}
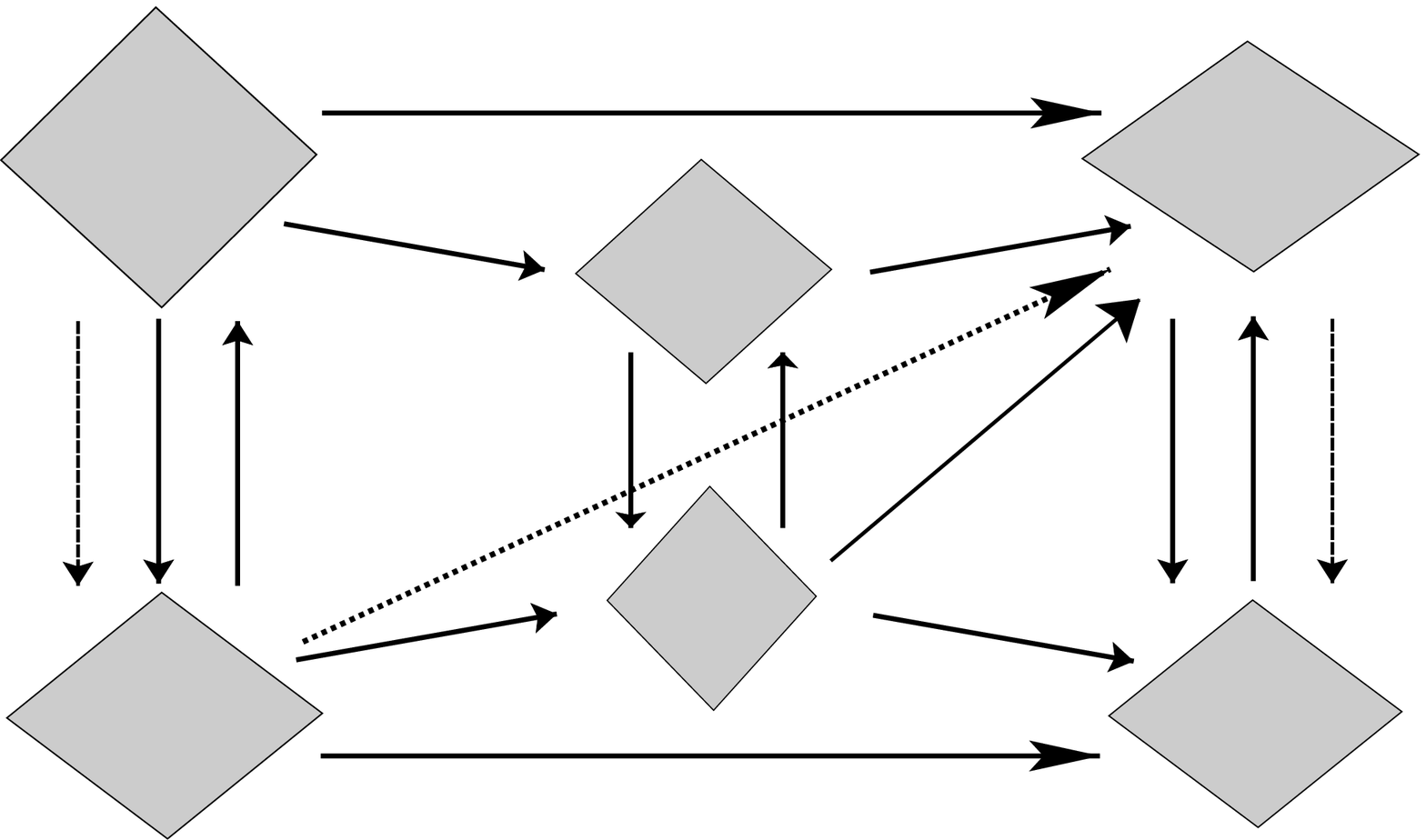
\end{center}
\caption{Definition of metric and stress tensors for finite growth.  \emph{Current metric} in Lagrangian configuration $\BC = \BF^T\Bg\BF$. \emph{Reference metric} in Eulerian configuration $\Bc = \BF^{-T}\BG\BF^{-1}$. The relationship between the Lagrangian-intermediate configuration and intermediate-Eulerian configuration can be defined $\BC^{e} = \BF^{g-T}\BC\BF^{g-1}$ and $\Btau = \BF^{e}\BS^{e}\BF^{eT}$ or directly $\Btau = \BF\BS\BF^T$, respectively \cite{Goktepe2010}.}
\label{fig3}
\end{figure}

The condition $J:=\det [\BF]>0$,  $J^e:=\det [\BF^e]>0$ and $J^g:=\det [\BF^g]>0$ provides the non-penetrable deformations. Furthermore, the reference $\SB_0$ and the spatial $\SB$ manifolds are locally furnished with the covariant reference $\BG$ and current $\Bg$ metric tensors in the neighborhoods $\eulN_{X}$ of $\BX$ and $\eulN_{x}$ of $\Bx$, respectively. For the mapping between the co- and contravariant objects in Lagrangian and Eulerian states, these metric tensors are needed \cite{marsden1983}. In addition, these metric tensors also show forth the map between the reference $\SB_0$ and the intermediate configuration $\SB_g$ and between the intermediate configuration $\SB_g$ and  the spatial state $\SB$.
Here we designate the \emph{right Cauchy Green tensor}, the inverse of the \emph{left Cauchy Green tensor}, and the \emph{elastic right Cauchy Green tensor},
\eb
\BC = \BF^T\Bg\BF,
\qquad \Bc = \BF^{-T} \BG \BF^{-1}  \qquad \text{and} \qquad
\BC^e = \BF^{eT}\Bg\BF^e,
\ee
as the covariant pull back of the spatial metric $\Bg$ and push forward of the Lagrangian metric $\BG$, respectively. Here, in the case of Cartesian coordinates, $\Bg$ and $\BG$ are the same as the identity tensor \cite{Goktepe2010}.
The \emph{left Cauchy Green tensor} or the \emph{Finger tensor} is defined by $\Bb=\Bc^{-1}$ .
For a geometric representation, we refer to Figure \ref{fig3}.
We define the Lagrangian unit vector $\Bn_0$ standing for the fiber direction in continuum as

\eb 
|\Bn_{0}|_{\BG}=1 \qquad\text{where}\qquad |\Bn_{0}|_{\BG}=\left(\Bn_{0}\cdot\BG\Bn_{0}\right )^{1/2}~.
\ee
Then through the tangent map deformed configuration can be formed as
\eb
\Bn=\BF^{e}\BF^{g}\Bn_0 ~.
\ee
The boundaries of the domain can be split into Dirichlet and Neumann type as $\partial \SB = \partial\SB^\varphi \cup \partial\SB^t$ and  $\partial\SB^{\varphi} \cap \partial\SB^t = \emptyset$.
\subsection{Governing equations for finite growth}\label{sec02-2}
The second law of thermodynamics states a positive dissipation by stress power and the objective rate of free energy. The dissipation inequality can be applied to form thermodynamically consistent stress relations as
\eb\label{dissipation}
\SD=\BS:\frac{1}{2}\dot\BC-\dot\psi = \hat\BP:\dot\BF -\dot\psi=\Btau:\frac{1}{2}\Lieder_{\Bv}\Bg-\dot\psi \ge 0 
\ee

It is significant to notice that  $\hat\BP$ implies the mixed-variant Piola tensor, which is represented by $\hat\BP:=\Bg\cdot\BP$ \cite{Goktepe2010}. The Lie derivative of the spatial metric, $\Lieder_{\Bv}\Bg=(\Bg\Bl+\Bl^T\Bg)$, is equivalent to the symmetric rate of deformation tensor. $\Bl=\dot\BF\BF^{-1}$ is the spatial velocity gradient. By the evaluation of the dissipation inequality, stress derivations can be performed. In finite growth, hyperelasticity constitutive equations are assumed to be held on the elastic state, and growth represents an intermediate stress-free state. Depending on the microstructure of the material or the tissue, the form of the growth tensor $\BF^{g}$ can be either isotropic or anisotropic. Stress arises in the solid domain due to the elastic part of the deformation gradient $\BF^{e}$. The growth does not have any energy contribution to the free energy function, meaning that purely growth tensor does not cause any stress in the body. The growth tensor can directly either be dependent on a scalar strain-driven variable or the microstructure that leads the growth evolution. It can be defined by evolution equations based on an internal variable. This type of growth depends on the growth criteria as it seen in the flow rule in plasticity formulation, and growth initiates when the mechanical driving force exceeds the threshold level. For the sake of convenience, in this study, it is assumed that every continuum point exceeds the threshold level in the solid body; hence growth tensor $\BF^{g}$ can be driven by certain scalar variables.
\begin{figure}[b!]
\centering
\def\svgwidth{0.7\textwidth}
\begingroup%
  \makeatletter%
  \providecommand\color[2][]{%
    \errmessage{(Inkscape) Color is used for the text in Inkscape, but the package 'color.sty' is not loaded}%
    \renewcommand\color[2][]{}%
  }%
  \providecommand\transparent[1]{%
    \errmessage{(Inkscape) Transparency is used (non-zero) for the text in Inkscape, but the package 'transparent.sty' is not loaded}%
    \renewcommand\transparent[1]{}%
  }%
  \providecommand\rotatebox[2]{#2}%
  \newcommand*\fsize{\dimexpr\f@size pt\relax}%
  \newcommand*\lineheight[1]{\fontsize{\fsize}{#1\fsize}\selectfont}%
  \ifx\svgwidth\undefined%
    \setlength{\unitlength}{764.29903766bp}%
    \ifx\svgscale\undefined%
      \relax%
    \else%
      \setlength{\unitlength}{\unitlength * \real{\svgscale}}%
    \fi%
  \else%
    \setlength{\unitlength}{\svgwidth}%
  \fi%
  \global\let\svgwidth\undefined%
  \global\let\svgscale\undefined%
  \makeatother%
  \begin{picture}(1,0.25346145)%
    \lineheight{1}%
    \setlength\tabcolsep{0pt}%
    \put(0,0){\includegraphics[width=\unitlength]{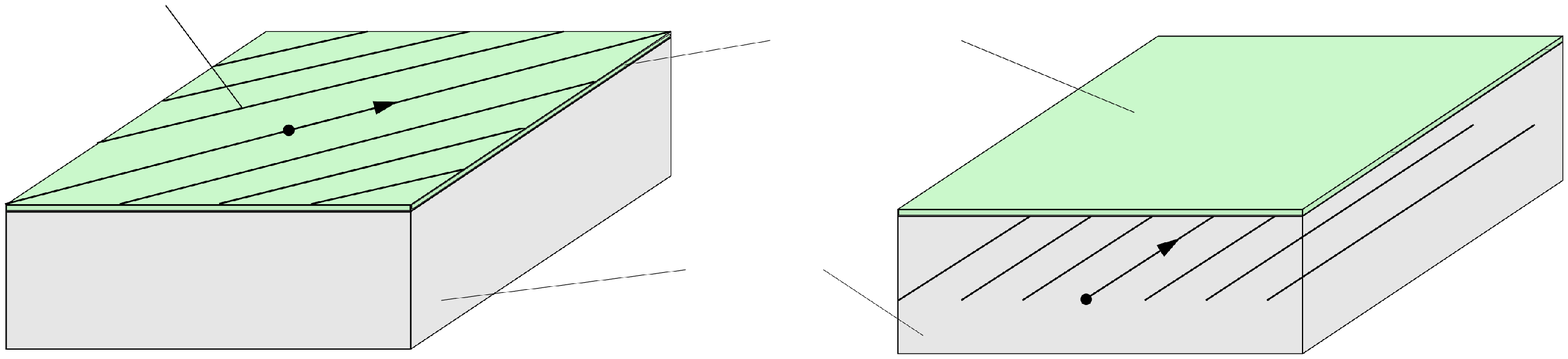}}%
    \put(0.50,0.21){\color[rgb]{0,0,0}\makebox(0,0)[lt]{\lineheight{1.25}\smash{\begin{tabular}[t]{l}  stiff film\end{tabular}}}}%
    \put(0.38,0.06){\color[rgb]{0,0,0}\makebox(0,0)[lt]{\lineheight{1.25}\smash{\begin{tabular}[t]{l} soft substrate\end{tabular}}}}%
    \put(0.02,0.24){\color[rgb]{0,0,0}\makebox(0,0)[lt]{\lineheight{1.25}\smash{\begin{tabular}[t]{l} stiff fiber\end{tabular}}}}%
    \put(0.16,0.16){\color[rgb]{0,0,0}\makebox(0,0)[lt]{\lineheight{1.25}\smash{\begin{tabular}[t]{l} $\Bn_0$ \end{tabular}}}}%
    \put(0.68,0.0){\color[rgb]{0,0,0}\makebox(0,0)[lt]{\lineheight{1.25}\smash{\begin{tabular}[t]{l} $\Bn_0$ \end{tabular}}}}%
  \end{picture}%
\endgroup%

\caption{Transverse isotropy:  stiff fibers embedded in the upper stiff film and stiff fibers embedded in the compliant  soft substrate, respectively}
\label{fig4}
\end{figure}
In general, fiber-reinforced rubber polymers, fibrous soft biological tissues, and reinforced composite elastomers exhibit nearly incompressible responses and, in the reinforcement direction, nearly inextensible behavior. It is common to observe that the growth of multi-layered materials (stiff bilayer film, hydrogels, etc.) or tissues (skin, artery, plant, etc.) can cause instabilities in the form of wrinkles, stripes, and secondary buckling shapes under residual stresses. The multi-layered materials can be composed of layers of different isotropic materials with different stiffnesses (e.g., stiff film on a compliant soft substrate). They can be formed by combinations of isotropic and fiber-reinforced anisotropic materials that we usually observe in nature, see  Figure \ref{fig4}. Depending on the physical conditions, these solids can be subjected to growth or shrinkage, resulting in compression or tension in the matrix and fibers, respectively. In line with the macro-level continuum approach, we introduce the elastic Helmholtz free energy function for one fiber family elastomers or soft tissues that is additively split into volumetric, isotropic, and anisotropic parts.
\eb \label{free-decomposed}
\psi(\Bg;\BF^e,\Bn_0):= \psi_{vol}(J^e) + \psi_{iso}(\Bg;\BF^e) + \psi_{ani}(\Bg;\BF^e,\Bn_0)
\ee
In general, stored free energy for hyperelastic solid is governed by three invariants of the right Cauchy Green tensor. These invariants can be used in combination with any kind of isotropic hyperelastic material model. In growth-induced hyperelasticity free energy formulation, elastic invariants are determined as
\eb\label{inv-iso}
I_1^e:= \tr\BC^e , \quad   I_2^e:=\frac{1}{2}\left[I_1^{e2}-\tr(\BC^{e2})\right],\quad \textrm{and}  \quad I_3^e:=\det \BC^e=J^{e2}
\ee
In the present study, an incompressible Neo-Hookean model, which considers only the first invariant for the isotropic response, is used to keep the formulation simple. The elastic right Cauchy Green tensor $\BC^{e} = \BF^{g-T}\BC\BF^{g-1}$ is defined in the intermediate configuration. In addition to the isotropic response, the anisotropic response that is generated by fibers requires additional invariants  $I_4$ and $I_5$. These additional invariants are introduced in terms of $\Bn_0$ unit vector of fibers in the undeformed configuration  
\eb\label{inv-ani}
I_4^e:= \Bn_0\cdot \BC^e \Bn_0 \, \qquad I_5^e=\Bn_0\cdot \BC^{e2} \Bn_0 \,,
\ee
that relates the energy storage due to fiber reinforcement in the material. The volumetric part of the free energy function in equation (\ref{free-decomposed}) is defined as
\eb\label{free-vol}
\psi_{vol}^e(J)= \frac{\lambda}{2}\left ( J^{e}-1 \right )^2
\ee
imposing the quasi-incompressibility, where ${\lambda}$ is Lame constant of the material. The second term of the free energy function is the isotropic part and represented by the Neo Hookean model as
\eb\label{free-iso}
\psi_{iso}^e(\Bg;\BF^e) = \frac{\mu_0}{2}(I_1^e-{2}\mbox{ln}\left (J^{e} \right )-3) \,,
\ee
where $\mu_0$ is the shear modulus. The last term of the free energy function stands for the anisotropic part
\eb\label{free-ani}
\psi_{ani}^e(\Bg;\BF^e,\Bn_0)=  \mu_f (I_4^e-1)^2
\ee
which is considered as a function of the fourth invariant. The volumetric and anisotropic parts of the free energy function will enter the formulation as constraints that lead the quasi-incompressible and quasi-inextensible behavior. In this study, fibers are assumed to carry both tension and compression.  In some applications, such as some soft biological tissues, it is assumed that fibers contribute only in tension but do not take any load in compression. For such applications, the anisotropic part of the free energy can be defined by Macaulay brackets $\langle\bullet\rangle$ that represent the tension-only behavior. Then, the parenthesis term in equation (\ref{free-ani}) can be replaced by $\langle I_4^e-1 \rangle^2$ for tension-only actions. By evaluating the dissipation inequality (\ref{dissipation}), the second Piola-Kirchhoff stress $\BS$, the Piola stress $\BP$ and the Kirchhoff stress $\Btau$, can be obtained thermodynamically conjugate to the right Cauchy Green deformation tensor $\BC$, and deformation gradient $\BF$, and the current metric $\Bg$, respectively \cite{Goktepe2010, Javili2015}. Then Kirchhoff stress $\Btau$ can be determined by the push forward of the second Piola-Kirchhoff stress $\BS$ as
\eb\label{overallstress}
\begin{array}{l@{}l}
\BS= 2 \p{\BC}{\psi} =2 \p{\BC^e}{\psi}:\p{\BC}{\BC^e} = \BF^{g-1}\cdot\BS^e\cdot\BF^{g-T}    \,, \\[2.5Ex]
\BP= \p{\BF}{\psi} =\p{\BF^e}{\psi}:\p{\BF}{\BF^e} = \BP^{e}\cdot\BF^{g-T}    \,, \\[2.5Ex]
\Btau= 2 \p{\Bg}{\psi} =\BF^{e}\cdot\BS^e\cdot\BF^{eT} 
\end{array}
\ee
Corresponding tangent moduli $\IC$, $\IA$, and $\IIc$ can be derived by push forward and pull back operations from elastic moduli $\IC^e$, $\IA^e$, and ${\IIc^e}$, which can be determined by taking the second derivative with respect to conjugate tensors $\BC$, $\BF$ and $\Bg$, respectively \cite{Goktepe2010,Javili2015}.
\eb\label{overallmoduli}
\begin{array}{l@{}l}
\IC= 2 \p{\BC}{\BS} =2 \p{\BC}( \BF^{g-1}\cdot\BS^e\cdot\BF^{g-T})= \left[ \BF^{g-1}\bar\dyad \BF^{g-1}  \right ]    : \IC^e : \left[ \BF^{g-T}\bar\dyad \BF^{g-T}  \right ]     \,, \\[2.5Ex]
\IA= \p{\BF}{\BP} =\p{\BF}( \BP^{e}\cdot\BF^{g-T})= \left[ \Bnone\bar\dyad \BF^{g-T}  \right ]    : \IA^e : \left[ \Bnone\bar\dyad \BF^{g-T}  \right ]  \,, \\[2.5Ex]
\IIc= 2\p{\Bg}{\Btau} = \left[ \BF\bar\dyad \BF  \right ]    : \IC : \left[ \BF^T\bar\dyad \BF^{T}  \right ]     \
\end{array}
\ee
where $\bar\dyad$ implies non-standard tensor product and can be defined as $\left[\bullet\bar\dyad\circ \right]_{ijkl} = \left[\bullet\right]_{ik} \left[\circ\right]_{jl}$. After defining Lagrangian, Eulerian, and two-point stress relations with elastic components in equation (\ref{overallstress}), it is required to derive corresponding elastic second Piola-Kirchhoff stress $\BS^e$ and the Kirchhoff stress $\Btau$. It should be noted that, elastic Kirchhoff stress can be determined by the push forward of the second Piola-Kirchhoff stress. However, as the outcome of the Eulerian configuration is more compact, we prefer to give both deformed and undeformed configurations in this study.
\subsection{Definition of the growth tensor}\label{sec02-3}

Depending on the nature of the material or the tissue, the growth tensor $\BF^g$ can be isotropic, transversely isotropic, orthotropic or anisotropic \cite{Kuhl2014}. The simplest approach is to express the growth tensor as isotropic; by this approach, the growth amount is equal in all directions. It can be a function of one scalar growth parameter \textit{g}, and then the isotropic growth tensor is defined as
\eb\label{isogrowth}
\BF^g = \left[1 +\textit{g} \right]\Bnone
\ee
where $\Bnone$ is the identity tensor. Note that, \textit{g} is the scalar growth parameter here and should not be confused by metric tensor $\Bg$. If the growth parameter g is zero, the growth tensor becomes equal to identity ($\BF^g=\Bnone$) which means there is no growth. Parameter \textit{g} can be positive or negative, corresponding to growth or shrinkage in the solid, respectively.
In this study, we consider a special case of growth, namely transversely isotropic or planar growth. It exhibits uneven growth for different orthogonal axes. Here, we focus on growth in the membrane plane while there is no growth in the thickness direction. In this way, the growth tensor $\BF^g$ is defined as
\eb\label{anisogrowth}
\BF^g = \left[1 +\textit{g} \right]\Bnone - \textit{g}\left[\Bm_0\dyad\Bm_0 \right]
\ee
where $\Bm_0$ is the unit normal of membrane in the reference configuration. Equation (\ref{anisogrowth}) describes an in plane growth within membrane and there is no growth along $\Bm_0$ direction. For the transversely isotropic growth, the inverse of the growth tensor $\BF^{g-1}$ can be defined as below \cite{Raush2014}
\eb\label{invanisogrowth}
\BF^{g-1} = \left[\frac{1}{1+\textit{g}}\right]\Bnone - \frac{1}{\textit{g}}\left[\Bm_0\dyad\Bm_0 \right]
\ee
\subsection{Stresses and moduli expressions in Lagrangian configuration}\label{sec02-4}
Similar to the free energy decomposition in (\ref{free-decomposed}), the elastic second Piola-Kirchhoff stress $\BS^{e}$ is defined as additively decoupled terms, namely volumetric, isochoric, and anisotropic such as
\eb\label{elastic2piokirchhoff}
\BS^e= 2 \p{\BC^e}{\psi^e} = \BS^{e}_{vol} + \BS^{e}_{iso} + \BS^{e}_{ani}
\ee
In line with the (\ref{elastic2piokirchhoff}), the elastic Lagrangian moduli $\IC^{e}$, can be additively split into volumetric, isotropic, and anisotropic parts respectively.
\eb\label{elasticLagMod}
\IC^e= 2 \p{\BC^e}{\BS^e} = 4 \p{\BC^e\BC^e}{\psi^e} = \IC^{e}_{vol}  + \IC^{e}_{iso} + \IC^{e}_{ani} 
\ee
The volumetric part of the elastic second Piola-Kirchhoff stress $\BS^{e \textit{ } vol}$ reads
\eb
\BS^{e}_{vol} :=2\p{\BC^e}\psi^{e}_{vol}(J^e)=p^e\BC^{e-1}\quad \text{with}\quad 
p^e:=J^e\psi'^{e \textit{ } vol}(J^e) = J^e{\lambda}(J^e-1)\,
\ee
The Lagrangian moduli expression for the volumetric part is determined as
\eb
\begin{array}{l@{}l}
\IC^{e}_{vol} :=2\p{\BC^{e}_{vol}}\BS^{e}_{vol}(J^e)=(p^e+\hat\kappa_L)\BC^{e-1}\dyad
\BC^{e-1} - {2}p^e\II_{\BC^{e-1}} \,, \\[2.5Ex]
\text{with}\quad 
\hat \kappa_L =J^{e \textit{ } 2}\psi''^{e}_{vol}(J^e)={\lambda}J^{e \textit{ } 2}\,
\end{array}
\ee
where $\II_{\BC^{e-1}}^{ABCD}= \frac{1}{2}\left[\BC_{AC}^{e-1}\BC_{BD}^{e-1}+\BC_{AD}^{e-1}\BC_{BC}^{e-1} \right]$. The isotropic elastic second Piola-Kirchhoff stress is obtained as
\eb
\begin{array}{l@{}l}
\BS^{e}_{iso}:=2\p{\BC^e}\psi^{e}_{iso}(\BC^e)= \mu_0(\BG^{-1}-\BC^{e-1})
\end{array}
\ee
Then, the Lagrangian elastic moduli for the isotropic response can be derived as
\eb
\begin{array}{l@{}l}
\IC^{e}_{iso}:=2\p{\BC^{e}}\BS^{e}_{iso}(\BC^e) = 2\mu_{0}\II_{\BC^{e-1}}
\end{array}
\ee
The last term of the elastic second Piola-Kirchhoff stress is the anisotropic response due to fiber reinforcement that can be derived as
\eb
\BS^{e}_{ani}:=2\p{\BC^e}\psi^{e}_{ani}(I_4)=4\mu_f \left(I_4^e-1\right)\Bn_0\dyad\Bn_0
\ee
Finally, the corresponding anisotropic elastic Lagrangian moduli read
\eb
\IC^{e}_{ani}:=2\p{\BC^e}\BS^{e}_{ani}(I_4)=8\mu_f\Bn_0\dyad\Bn_0\dyad\Bn_0\dyad\Bn_0
\ee
\subsection{Stresses and moduli expressions in Eulerian configuration}\label{sec02-5}

Similar to Lagrangian configuration, the Kirchhoff stress can be defined by additively decomposition of volumetric, isotropic, and anisotropic terms
\eb
\Btau:=2\p{\Bg}{\psi}=\Btau_{vol}+\Btau_{iso}+\Btau_{ani}
\ee
and the corresponding Eulerian tangent moduli read
\eb
\IIc:=4\p{\Bg\Bg}^2 \psi(\Bg;\BF,\Bf_0)=\IIc_{vol}+\IIc_{iso}+\IIc_{ani}
\ee
The Kirchoff stress for the volumetric response can be derived as
\eb
\Btau_{vol}:=2\p{\Bg}\psi_{vol}(J)=p\Bg^{-1}\quad \text{with}\quad 
p:=J\psi'_{vol}(J) = J{\lambda}(J-1)\,.
\ee
The Eulerian tangent moduli expression for the volumetric part is defined as
\eb
\IIc_{vol}:=4\p{\Bg\Bg}^2 U(J)=(p+\hat\kappa)\Bg^{-1}\dyad \Bg^{-1}-2p\II\,,\quad \text{with} \quad \hat \kappa =J^2\psi''_{vol}(J)={\lambda}(J^2+1)\,.
\ee
The elastic Kirchhoff stress for the isotropic response reads
\eb \label{def-tau-iso}
\Btau_{iso}:=2\p{\Bg}{\psi_{iso}(\Bg;\BF)}=\mu_0(\Bb-\Bg^{-1})
\ee
Then the corresponding elastic Eulerian moduli for the isotropic response can be derived as
\eb
\IIc_{iso}:=4\p{\Bg\Bg}^2\psi_{iso}(\Bg;\BF)
=2\mu_0\II
\ee
Anisotropic part of the Kirchhoff stress representing the fiber reinforcement is defined as
\eb
\Btau_{ani}=2\p{\Bg}\psi_{ani}(\Bg;\BF,\Bn_0)=4\mu_f( I_4-1)\Bn\otimes \Bn \,.
\ee
Moreover, the corresponding anisotropic part of Eulerian moduli reads
\eb
\IIc_{ani}=4\p{\Bg\Bg}^2\psi_{ani}(\Bg;\BF,\Bn_0)=8\mu_f\Bn\otimes\Bn\otimes\Bn\otimes\Bn\,.
\ee
Finally, by using the relations of pull back and push forward operations defined in  (\ref{overallstress}) and (\ref{overallmoduli}), one can end up with the overall stress and moduli in terms of the second Piola-Kirchhoff stress $\BS$, the Piola stress $\BP$, the Kirchhoff stress $\Btau$ and the corresponding moduli $\IC$, $\IA$ and ${\IIc}$, respectively.

\section{Variational formulation for anisotropic and incompressible continuum } \label{sec03}
\subsection{Variational formulation for finite elasticity}\label{sec03-1}
The potential functional can be defined as below for the \emph{finite elasticity}
\eb
\hat\Pi(\Bvarphi,t) := \hat\Pi^{int}(\Bvarphi,t)-  \hat\Pi^{ext}(\Bvarphi,t)\,,
\label{var1}
\ee 
where 
\eb
\hat\Pi^{int}(\Bvarphi,t) := \nint{}{\SB}\psi(\Bg,\BF^{e})\dV  \AND \hat\Pi^{ext}(\Bvarphi):=  \nint{}{\SB}\Bvarphi\cdot\rho_{0}\bar\Bgamma  \dV
+ \nint{}{\partial\SB_t}\Bvarphi\cdot\bar\BT \dA\,.
\ee
For elastic loading, the stored energy in the body can be defined by $\hat\Pi^{int}(\Bvarphi)$, and $\hat\Pi^{ext}(\Bvarphi)$ refers to the work done by external forces. In equation (\ref{var1}),  $\rho_0$, $\bar\Bgamma$, and $\bar\BT$ are  the density, body force, and the surface traction, respectively.  $\psi(\Bg,\BF^{e})$ is the volume-specific elastic Helmholtz free energy. The boundary value problem governing finite elasticity is obtained from the elastic potential by the \emph{principle of minimum potential energy} in the variational form
\begin{equation}\label{inf1}
\Bvarphi_t = \Arg \{\inf_{\Bvarphi_t\in\SW} \hat\Pi(\Bvarphi,t)\}
\end{equation}
subject to Dirichlet-type boundary condition
\eb\label{dirichlet}
\SW = \{ \Bvarphi_t ~|~ \Bvarphi_t \in \SB \quad \land \quad \Bvarphi_t = \bar \Bvarphi \quad\textrm{on}\quad \partial\SB_u  \}\,.
\ee
Due to the stationary condition of the elastic potential $\hat\Pi(\Bvarphi,t)$, the variation of (\ref{inf1}) along with localization theorem yields the \emph{Euler--Lagrange equation}
\eb \label{bom}
J^{e}\div [J^{e-1}\Btau] +\rho_0 \bar\Bgamma  = 0
\ee
yielding to  the balance of linear momentum for static problems in the domain $\SB$  along with Neumann-type boundary condition
\eb
\BP\cdot\BN=\Btau\cdot \Bn = \bar\BT \qquad \textrm{on} \qquad  \partial\SB_t
\ee
where we have used the identity of \emph{Nanson's formula} as $J\BF^{-T}\BN~dA=\Bn~da$. The element formulation can be derived by the consistent linearization of the weak form obtained as the first variation of (\ref{inf1}).
\subsection{A mixed variational formulation for quasi-incompressible and quasi-inextensible continuum}\label{sec03-2}
The quasi-incompressible and quasi-inextensible behavior can be enforced by two additional penalty terms in the decomposed representation \req{free-decomposed} of the free energy functional as
\eb \label{var2}
\hat\Pi(\Bvarphi, p^{e}, \theta,s^{e},\lambda) :=   
\nint{}{\SB}\pi^{\ast}(\Bvarphi, p^{e}, \theta,s^{e},\lambda)\dV
-\hat\Pi^{ext}(\Bvarphi,t) \,.
\ee
The mixed potential density in equation \req{var2} is defined as
\eb\label{pot-den} 
\pi_{int}^{\ast}(\Bvarphi, p^{e}, \theta,s^{e},\lambda)=
\psi_{iso}(\Bg,\BF^{e})+ \underbrace{p^{e}(J^{e}-\theta)+\psi_{vol}(\theta)}_{\textrm{volumetric constraint}}
+\underbrace{s^{e}(I_4^{e}-\lambda)+\psi_{ani}(\lambda)}_{\textrm{inextensibility constraint}}\,.
\ee
Here,  $p^{e},\,s^{e}$ are penalty parameters dual to the kinematic quantities $\theta,\,\lambda$. The deformation of the body enforced by incompressible and inextensible constraints is governed by the \emph{mixed saddle point  principle}
\begin{equation}\label{inf2}
\{\Bvarphi_t,\theta,p^{e},\lambda,s^{e} \} = \Arg\{\inf_{\Bvarphi_t  \in \SW}\inf_{\theta}\inf_{\lambda}\sup_{p^{e}}\sup_{s^{e}} \hat\Pi(\Bvarphi,t)\} \,.
\end{equation}
subject to the boundary conditions $\SW = \{ \Bvarphi_t ~|~ \Bvarphi_t \in \SB \quad \land \quad \Bvarphi_t = \bar \Bvarphi \quad\textrm{on}\quad \partial\SB_u  \}$.
Taking the first variation of  (\ref{var2}) with respect to $\Bvarphi$, $p^{e}$, $\theta$, $s^{e}$ and $\lambda$ provides the weak form
\eb\label{weak1}
\begin{array}{l@{~=~}l}
\delta_\Bvarphi \hat\Pi(\Bvarphi, p^{e}, \theta,s^{e},\lambda) & \nint{}{\SB}\{\left(\Btau_{iso}+ p^{e} J^{e} \Bg^{-1} + 2s^{e} \Bn\dyad\Bn \right):
\frac{1}{2}\Lieder_{\delta\Bvarphi}\Bg\}\dV	- \delta\hat\Pi^{ext}(\Bvarphi) = 0\,, \\[2.5Ex]
\delta_{p^{e}} \hat\Pi(\Bvarphi, p^{e}, \theta,s^{e},\lambda) & \nint{}{\SB} \delta p^{e}~ (J^{e}-\theta) \dV = 0\,, \\[2.5Ex]
\delta_\theta \hat\Pi(\Bvarphi, p^{e}, \theta,s^{e},\lambda) & \nint{}{\SB} \delta \theta~ (\psi_{vol}'^{e}(\theta)- p^{e} ) \dV = 0\,, \\[2.5Ex]
\delta_{s^{e}} \hat\Pi(\Bvarphi, p^{e}, \theta,s^{e},\lambda) & \nint{}{\SB} \delta s^{e}~ (I_4^{e}-\lambda) \dV = 0\,, \\[2.5Ex]
\delta_ \lambda \hat\Pi(\Bvarphi, p^{e}, \theta,s^{e},\lambda) & \nint{}{\SB} \delta \lambda ~(\psi_{ani}'^{e}(\lambda)- s^{e} ) \dV = 0\,, 
\end{array}
\ee
from which the mixed finite element method can be formulated. Therein, $\Lieder_{\delta\Bvarphi}\Bg$ is the  Lie derivative of the current metric along the variation $\delta\Bvarphi$.
Taking the variation of the potential density \req{pot-den}, \emph{Euler--Lagrange equations} of the mixed variational principle gives
\eb\label{eul-lag}
\begin{array}{ll@{~=~}l}
1.\qquad & J\div [J^{-1}\Btau] +\rho_0 \bar\Bgamma & 0 \\[1Ex]
2. & J - \theta & 0 \\[1Ex]
3. & \psi_{vol}'(\theta)- p^{e} & 0 \\[1Ex]
4. &  I_4-\lambda& 0 \\[1Ex]
5. &  \psi_{ani}'(\lambda)- s^{e} & 0 \\[1Ex]
\end{array}
\ee
along with the Neumann-type boundary conditions
$\SW_t=\{\Bsigma\cdot\Bn=\Bt \quad \textrm{on} \quad \partial\SB_t\}$. The consistent linearization of the mixed potential has been implemented and discretized in Python-based open-source finite element platform \emph{FEniCS} \cite{Logg2012}. The linear Newton iterations are solved through MUltifrontal Massively Parallel sparse direct Solver (MUMPS)\cite{Amestoy2000}.
\section{Representative numerical examples}\label{sec04}
A numerical example is demonstrated in this section through a boundary value problem. As a mathematical representative case, effects of fiber's stiffness to the buckling mode under planar growth are studied using \emph{T2P0F0} element in the open-source automated finite element program  \emph{FEniCS} for a 3D stiff film on a compliant substrate.
\subsection{Planar growth of fiber reinforced 3D-bilayer stiff film on a compliant substrate}
\begin{figure}[thb!]
\centering
\def\svgwidth{0.6\textwidth}
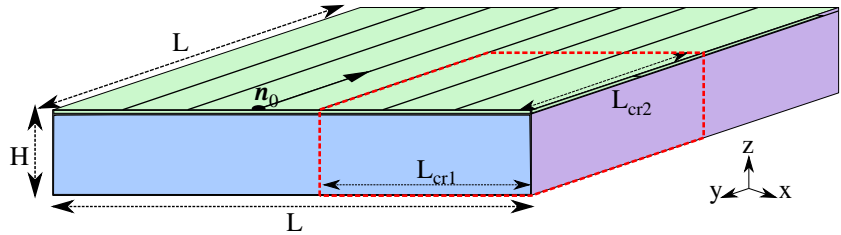
\caption{Schematic representation of the 3D bilayer plate (stiff transversely isotropic film on isotropic compliant substrate) with its geometric dimensions. L is the edge length of the square in-plane section, H is the total height including the film and the substrate, $\Bn_{0}$ is the stiff fiber direction along y axis, L\textsubscript{cr1} and L\textsubscript{cr2} are the critical lowest dimensions needed to be taken into account to reflect periodic behaviour of RVE (Representative Volume Element)}
\label{fig5}
\end{figure}
This example demonstrates the effect of the stiffness of the fibers, which exist in the stiff film, on the critical growth parameter of instability state under the planar growth. Since the analyzed domain is just a small part/patch of the real structure like a representative volume element, it is necessary to define boundary conditions properly to represent the overall behavior accurately. The schematic illustration of the 3D bilayer plate model is shown in Figure \ref{fig5}. In order to reflect the infinitely long plate analogy along x and y-axes, we assume that buckling the behavior of the structure repeats itself periodically with a critical wavelength ($\lambda_{cr}$).  Thus, we define periodic boundary conditions on sidewalls (blue and purple faces and opposite faces in \ref{fig5}).  To define periodic boundary conditions, we relate these face pairs by displacement relations on the left face at $\Bx=0$ and on the right face at $\Bx=L_{cr1}$ as shown below:
\eb\label{periodicBCx}
\begin{array}{l@{}l}
\Bx^L -  \Bx^R =L_{cr1} \\[2.5Ex]
\By^L =  \By^R \qquad \qquad \text{and} \qquad \Bu^L-\Bu^R=0 \\[2.5Ex]
\Bz^L =  \Bz^R 
\end{array}
\ee
Similary, we define periodic boundary relations on the front face at $\By=0$ and on the rear face at $\By=L_{cr2}$   as show below:
\eb\label{periodicBCy}
\begin{array}{l@{}l}
\Bx^F =  \Bx^B \\[2.5Ex]
\By^F -  \By^B =L_{cr2} \quad \quad \text{and} \qquad \Bu^F-\Bu^B=0\\[2.5Ex]
\Bz^F =  \Bz^B 
\end{array}
\ee
where superscripts $R$, $L$, $F$, and $B$ refer to left, right, front, and back sides, respectively. Since an infinitely long and wide plate is modeled as a finite-sized computational domain with periodic boundary conditions, the fiber direction is no longer significant for the one family fiber configuration. In other words, for single-family fiber models with periodic boundary conditions, the same buckling modes are obtained for different angles. Therefore, the initial fiber angle $\Bn_{0}$ is aligned to the y-direction as  $\Bn_{0} =[0 \quad 1\quad0]^{T}$. In addition to periodic boundary conditions on sides, the bottom face of the substrate is fixed in all directions.

\begin{figure}[thb!]
\centering
\def\svgwidth{1\textwidth}
\begingroup%
  \makeatletter%
  \providecommand\color[2][]{%
    \errmessage{(Inkscape) Color is used for the text in Inkscape, but the package 'color.sty' is not loaded}%
    \renewcommand\color[2][]{}%
  }%
  \providecommand\transparent[1]{%
    \errmessage{(Inkscape) Transparency is used (non-zero) for the text in Inkscape, but the package 'transparent.sty' is not loaded}%
    \renewcommand\transparent[1]{}%
  }%
  \providecommand\rotatebox[2]{#2}%
  \newcommand*\fsize{\dimexpr\f@size pt\relax}%
  \newcommand*\lineheight[1]{\fontsize{\fsize}{#1\fsize}\selectfont}%
  \ifx\svgwidth\undefined%
    \setlength{\unitlength}{2809.23420583bp}%
    \ifx\svgscale\undefined%
      \relax%
    \else%
      \setlength{\unitlength}{\unitlength * \real{\svgscale}}%
    \fi%
  \else%
    \setlength{\unitlength}{\svgwidth}%
  \fi%
  \global\let\svgwidth\undefined%
  \global\let\svgscale\undefined%
  \makeatother%
  \begin{picture}(1,0.21094362)%
    \lineheight{1}%
    \setlength\tabcolsep{0pt}%
    \put(0,0){\includegraphics[width=\unitlength]{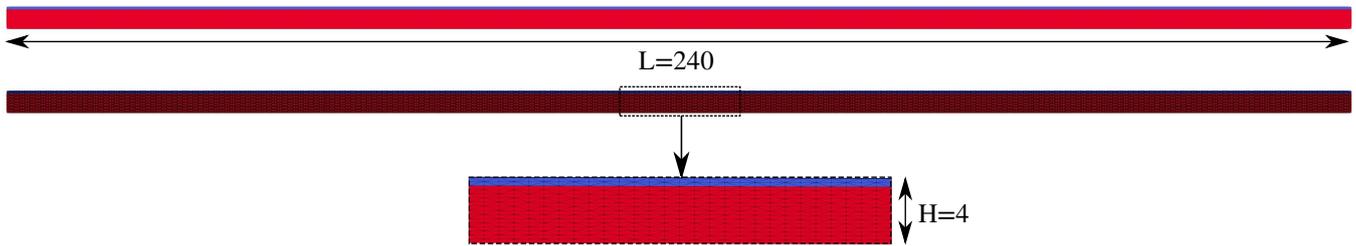}}%
    \put(0.46960464,0.13){\color[rgb]{0,0,0}\makebox(0,0)[lt]{\lineheight{1.25}\smash{\begin{tabular}[t]{l} L=240\end{tabular}}}}%
    \put(0.67580506,0.017){\color[rgb]{0,0,0}\makebox(0,0)[lt]{\lineheight{1.25}\smash{\begin{tabular}[t]{l} H=4\end{tabular}}}}%
  \end{picture}%
\endgroup%

\caption{Representation of the long (L=240 unit) and thin (W=1 unit) bilayer plate consist of substrate (red) and stiff film (blue), finite element mesh and zoomed view of FE mesh respectively.}
\label{fig6}
\end{figure}

In the definition of periodic boundary condition, it is essential to determine the critical (characteristic) geometric dimensions; it can also be introduced as the minimum size of RVE, which reflects the actual periodic behavior without enforcing any constraint to the overall structural stiffness. The minimum critical  lengths along the x and y-axes are shown in Figure \ref{fig5} as L\textsubscript{cr1} and L\textsubscript{cr2}. If the RVE size is not properly chosen, the computations could miss or artificially enforce a buckling mode. In order to capture the minimum required RVE size, three-dimensional long but thin bilayer plate models are defined along the fiber direction such that L=240 unit, H=4 unit, and W=1 unit dimensions, see Figure \ref{fig6}. The total height H =4 units, where 3.5 units correspond to the substrate and 0.5 units to the stiff film. The aim of the three-dimensional long but thin bilayer plate model is to identify characteristic wavelength ($\lambda_{cr}$) in a long regime by decreasing the effect of wall boundary conditions. Since single-family fiber reinforcement is studied in this study, fibers are defined along the y-axis. Then periodic boundary conditions for left-right and front-back face couples are defined as (\ref{periodicBCx}) and (\ref{periodicBCy}). In this study, the ratio of the shear modulus of the film layer ($\mu_{f}$)  to the substrate ($\mu_{s}$) is 100, which also provides a larger wavelength \cite{Dortdivanoglu2017}. 
Since this study aims to observe the effect of in-plane aligned fibers on bilayer three-dimensional plate buckling, fiber stiffness is the key parameter to be tested for a wide range. Accordingly, fiber stiffness values are taken as $\mu_{fiber}$={100, 250, 750, 1000, 2500. The three-dimensional bilayer plate is subjected to planar growth on both layers, and the growth parameter $g$ monotonically increases with time step. It is important to capture buckling modes of the bilayer plate with the proper time step increment $\Delta t$. The initial time step is defined as $1 \times 10^{-4}$, and it is divided by two when a convergence problem is encountered. After the first buckling initiation, the time step is also initialized to $1 \times 10^{-4}$ and it is kept constant until it reaches the secondary buckling stage. This process continues until the Newton-Raphson algorithm does not converge within 20 steps, even at the fifth loop of the time step division process. In order to trigger buckling at the critical growth value, a perturbation needs to be applied. In this study, the perturbation is defined as a minimal eccentric distributed load along x and y-directions. The material parameters used in the three-dimensional bilayer plate are given in Table \ref{mat-par}.

\begin{table}[thb!]
\centering
\caption{Material parameters used in the analysis of bilayer three-dimensional plate.} 
\label{mat-par}
\begin{tabular}{rllrll} 
\hline                             \\ [-2.5ex]
Parameter & Value & Unit & Parameter & Value & Unit  \\ \hline \\[-2Ex] 
 $\mu_{film}$  & $10^{2}$  & [~--~]&
$\mu_{subs}$ & \white{+}$1.0$  &[~--~]   \\ [0.5Ex]
$\kappa_{film}$  & $10^{5}$  & [~--~]&
$\kappa_{subs}$ & \white{+}$10^{3}$  &[~--~]   \\ [0.5Ex]
$\mu_{fiber}$   & ${100, 250, 750, 1000, 2500}$  &[~--~]   &
$\Bn_0$   &  \white{+}$[0,1,0]$  &[~--~]   \\ [0.5Ex]
\hline  
\end{tabular}
\end{table}

In order to determine the critical buckling wavelength ($\lambda_{cr}$) and related RVE dimensions, three-dimensional long but  thin bilayer plate geometry (see Figure \ref{fig6}) was analyzed for each fiber stiffness ($\mu_{fiber}$) values listed in Table \ref{mat-par}. It is aimed to observe the periodicity of the buckling behavior for all of the fiber stiffnesses and to identify the critical/characteristic lengths for the RVE that is needed for post-buckling analysis. The first buckling mode of bilayer plate without and with fibers has been observed and shown in Figure \ref{fig7}.
\begin{figure}[thb!]
\raggedleft
\vspace{7pt}%
\def\svgwidth{0.85\textwidth}
\begingroup%
  \makeatletter%
  \providecommand\color[2][]{%
    \errmessage{(Inkscape) Color is used for the text in Inkscape, but the package 'color.sty' is not loaded}%
    \renewcommand\color[2][]{}%
  }%
  \providecommand\transparent[1]{%
    \errmessage{(Inkscape) Transparency is used (non-zero) for the text in Inkscape, but the package 'transparent.sty' is not loaded}%
    \renewcommand\transparent[1]{}%
  }%
  \providecommand\rotatebox[2]{#2}%
  \newcommand*\fsize{\dimexpr\f@size pt\relax}%
  \newcommand*\lineheight[1]{\fontsize{\fsize}{#1\fsize}\selectfont}%
  \ifx\svgwidth\undefined%
    \setlength{\unitlength}{3129.20556256bp}%
    \ifx\svgscale\undefined%
      \relax%
    \else%
      \setlength{\unitlength}{\unitlength * \real{\svgscale}}%
    \fi%
  \else%
    \setlength{\unitlength}{\svgwidth}%
  \fi%
  \global\let\svgwidth\undefined%
  \global\let\svgscale\undefined%
  \makeatother%
  \begin{picture}(1,0.30312268)%
    \lineheight{1}%
    \setlength\tabcolsep{0pt}%
    \put(0,0){\includegraphics[width=\unitlength]{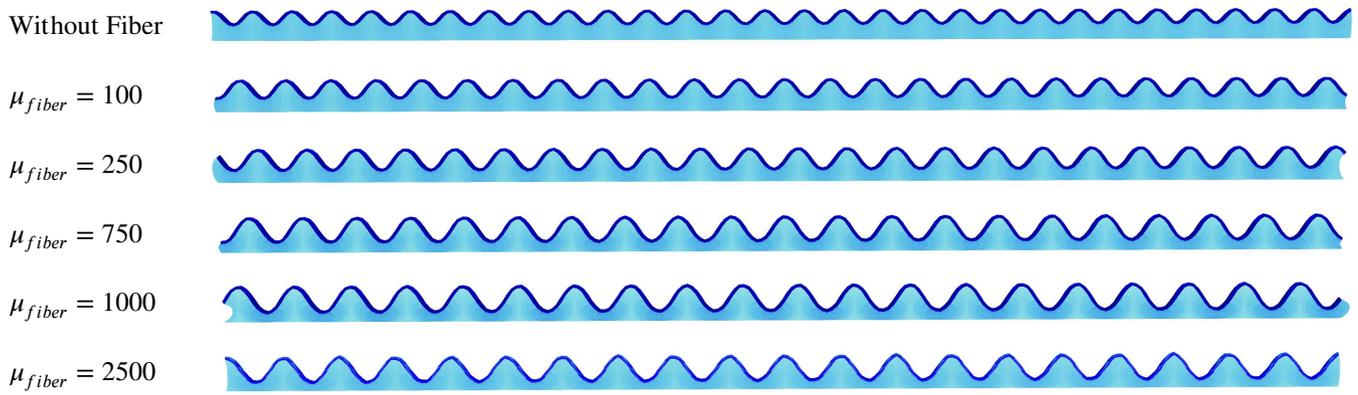}}%
    \put(-0.17,0.31){\color[rgb]{0,0,0}\makebox(0,0)[lt]{\lineheight{1.25}\smash{\begin{tabular}[t]{l}$\mathrm{Without \ Fiber}$\end{tabular}}}}%
    \put(-0.17,0.25){\color[rgb]{0,0,0}\makebox(0,0)[lt]{\lineheight{1.25}\smash{\begin{tabular}[t]{l}$\mu_{fiber}=100$\end{tabular}}}}%
    \put(-0.17,0.19){\color[rgb]{0,0,0}\makebox(0,0)[lt]{\lineheight{1.25}\smash{\begin{tabular}[t]{l}$\mu_{fiber}=250$\end{tabular}}}}%
    \put(-0.17,0.13){\color[rgb]{0,0,0}\makebox(0,0)[lt]{\lineheight{1.25}\smash{\begin{tabular}[t]{l}$\mu_{fiber}=750$\end{tabular}}}}%
    \put(-0.17,0.07){\color[rgb]{0,0,0}\makebox(0,0)[lt]{\lineheight{1.25}\smash{\begin{tabular}[t]{l}$\mu_{fiber}=1000$\end{tabular}}}}%
    \put(-0.17,0.01){\color[rgb]{0,0,0}\makebox(0,0)[lt]{\lineheight{1.25}\smash{\begin{tabular}[t]{l}$\mu_{fiber}=2500$\end{tabular}}}}%
  \end{picture}%
\endgroup%

\caption{First buckling mode shapes of 3D long and thin bilayer plates with different fiber stiffness reinforcement (scaled by x5).}
\label{fig7}
\end{figure}
Due to the contribution of fiber stiffness to the directional stiffness of the film, each model in Figure \ref{fig7} leads to a different critial growth value $g$ and a different wavelength. Figure \ref{fig7} shows the buckling mode shapes of long but thin bilayer plate for various fiber stiffness values. The wavelength information acquired from Figure \ref{fig7} will be used for the three-dimensional rectangular in-plane section of the bilayer plate as shown in Figure \ref{fig5} to examine the first buckling mode and post-buckling behavior. As shown in Figure \ref{fig7}, the number of wrinkles decreases by increasing the fiber stiffness $\mu_{fiber}$. It also means that the characteristic/critical wavelength increases with the fiber stiffness $\mu_{fiber}$. In order to reflect infinite plate behavior with periodic boundary conditions and capture at least a single wrinkle in the three-dimensional RVE, the minimum critical length L\textsubscript{cr2} should be at least  ($\lambda_{cr}$) unit. It can also be concluded that harmonics of ($\lambda_{cr}$) will also capture the exact buckling shape. Since the buckling behavior in the x-direction (perpendicular to fibers) is likely to be decoupled from fibers and represent a behavior similar to the no fiber case, it requires a smaller wavelength to capture post-buckling mode in the x-direction. Therefore, for the bilayer model shown in Figure \ref{fig5},  L\textsubscript{cr1} is taken as the same as  L\textsubscript{cr2} making the bilayer plate have a square in-plane section. To enlarge the coverage of the buckling more precisely within two wrinkles, L\textsubscript{cr} is set to 2$\lambda_{cr}$, which makes the geometrical dimensions 2$\lambda_{cr}$x2$\lambda_{cr}$xH as shown in Figure \ref{fig7b}.
\begin{figure}[thb!]
\centering
\vspace{15pt}%
\def\svgwidth{0.5\textwidth}
\begingroup%
  \makeatletter%
  \providecommand\color[2][]{%
    \errmessage{(Inkscape) Color is used for the text in Inkscape, but the package 'color.sty' is not loaded}%
    \renewcommand\color[2][]{}%
  }%
  \providecommand\transparent[1]{%
    \errmessage{(Inkscape) Transparency is used (non-zero) for the text in Inkscape, but the package 'transparent.sty' is not loaded}%
    \renewcommand\transparent[1]{}%
  }%
  \providecommand\rotatebox[2]{#2}%
  \newcommand*\fsize{\dimexpr\f@size pt\relax}%
  \newcommand*\lineheight[1]{\fontsize{\fsize}{#1\fsize}\selectfont}%
  \ifx\svgwidth\undefined%
    \setlength{\unitlength}{1625.41080985bp}%
    \ifx\svgscale\undefined%
      \relax%
    \else%
      \setlength{\unitlength}{\unitlength * \real{\svgscale}}%
    \fi%
  \else%
    \setlength{\unitlength}{\svgwidth}%
  \fi%
  \global\let\svgwidth\undefined%
  \global\let\svgscale\undefined%
  \makeatother%
  \begin{picture}(1,0.43505444)%
    \lineheight{1}%
    \setlength\tabcolsep{0pt}%
    \put(0,0){\includegraphics[width=\unitlength]{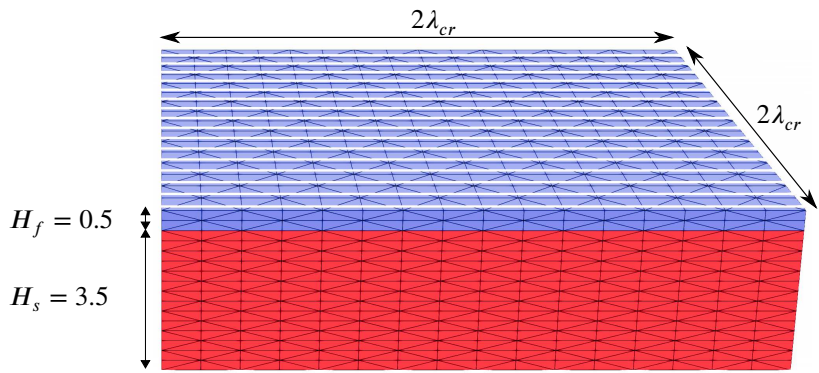}}%
    \put(0.4,0.52){\color[rgb]{0,0,0}\makebox(0,0)[lt]{\lineheight{1.25}\smash{\begin{tabular}[t]{l}$2\lambda_{cr}$\end{tabular}}}}%
    \put(0.91,0.38){\color[rgb]{0,0,0}\makebox(0,0)[lt]{\lineheight{1.25}\smash{\begin{tabular}[t]{l}$2\lambda_{cr}$\end{tabular}}}}%
    \put(-0.19365205,0.22616726){\color[rgb]{0,0,0}\makebox(0,0)[lt]{\lineheight{1.25}\smash{\begin{tabular}[t]{l}$H_{f}=0.5$\end{tabular}}}}%
    \put(-0.19365205,0.11301776){\color[rgb]{0,0,0}\makebox(0,0)[lt]{\lineheight{1.25}\smash{\begin{tabular}[t]{l}$H_{s}=3.5$\end{tabular}}}}%
  \end{picture}%
\endgroup%

\caption{Representation of three-dimensional RVE. White lines demonstrate the stiff fibers those are aligned to the y-direction}
\label{fig7b}
\end{figure}
In the light of these geometrical characteristic information, the effects of fiber stiffness over the primary buckling and the post-buckling regime was examined for the model shown in Figure \ref{fig7b}. Figure \ref{Fig9} shows pressure contours in a growth timeline where the structure is subjected to planar growth. It gives the results for isotropic bilayer plate and fiber-reinforced bilayer plate that have $\mu_{fiber}$=100, 250, 750, 1000 and 2500. It can be seen that in the isotropic case without fiber contribution, the first buckling triggered simultaneously in x and y-directions at the same time due to isotropic nature. The isotropic case yields to a labyrinth shape by increasing the growth parameter further, as reported in \cite{Huang2005}. For the primary buckling patterns of fiber-reinforced plates, it is clearly observed that buckling behaviors are similar, but the critical growth value that causes instability differs from each other. Furthermore, in the first buckling shape, sinusoidal wrinkles are initiated along the fiber direction, which causes energy relaxation of fibers in the critical growth parameter. Post-buckling behavior results as a secondary wrinkle with a different form and amplitude in the transverse direction perpendicular to fibers. Moreoever, it is also obtained that the form of wrinkles that are observed in the first buckling mode shifts from sinusoidal to triangular shape by increasing fiber stiffness. It can be noted from Figure \ref{Fig9} that during the buckling and post-buckling stage, while concave patterns have positive pressure contour, convex regions have negative pressure.  
\begin{figure}[h] 
\raggedright
\rotatebox{90} {$\mathrm{(a) Without \ fiber}$}
\def\svgwidth{0.90\textwidth}
\centering
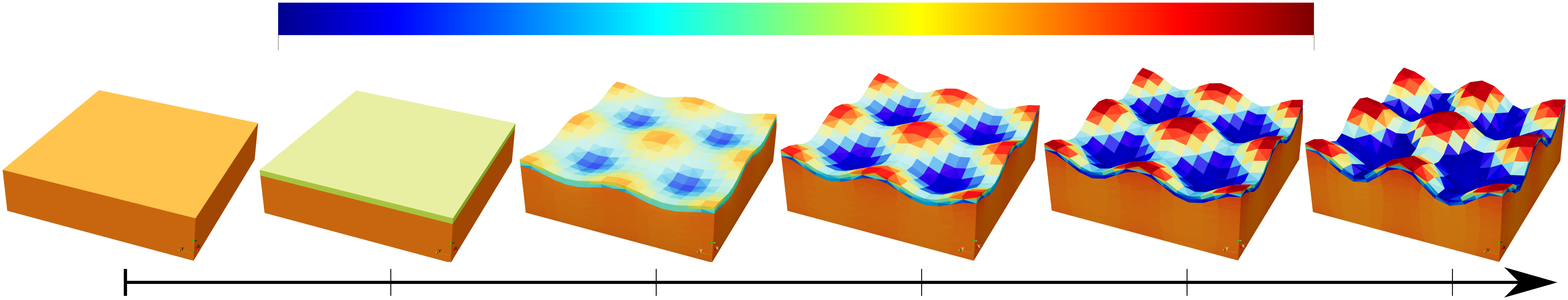
\label{Fig8a_NoFiber}
\end{figure}
\begin{figure}[h]
\raggedright
\rotatebox{90} {(b) $\mu_{fiber}=100$}
\def\svgwidth{0.9\textwidth}
\centering
\includegraphics*[width=0.90\textwidth]{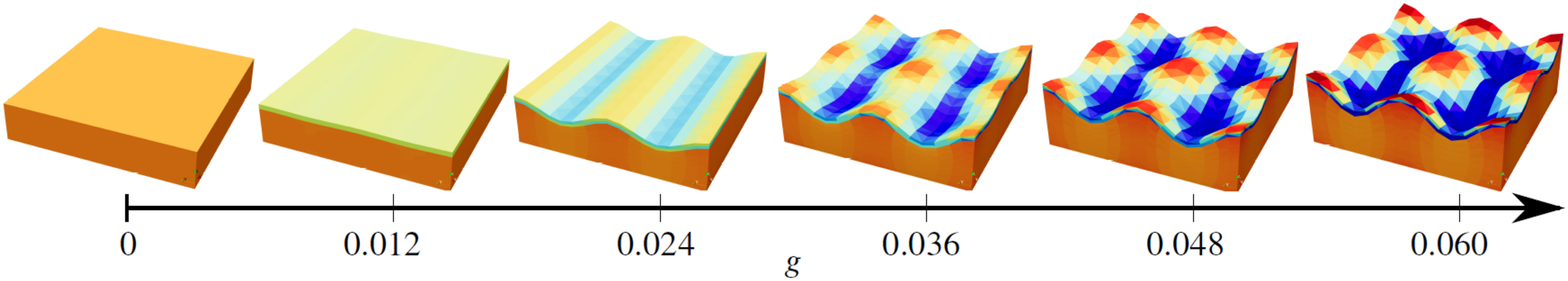}
\label{Fig8a_mu100}
\end{figure}
\begin{figure}[h]
\raggedright
\rotatebox{90} {(c) $\mu_{fiber}=250$}
\def\svgwidth{0.9\textwidth}
\centering
\includegraphics*[width=0.90\textwidth]{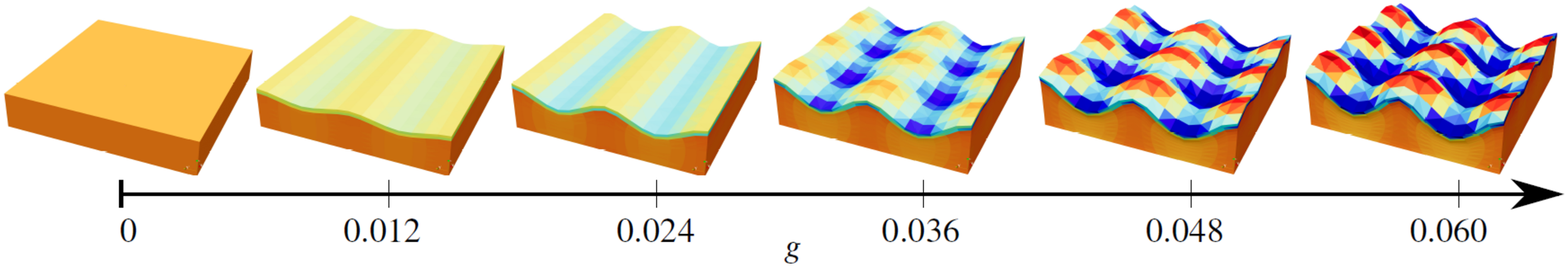}
\label{Fig8a_mu250}
\end{figure}
\begin{figure}[!h]
\raggedright
\rotatebox{90} {(d) $\mu_{fiber}=750$}
\def\svgwidth{0.9\textwidth}
\centering
\includegraphics*[width=0.90\textwidth]{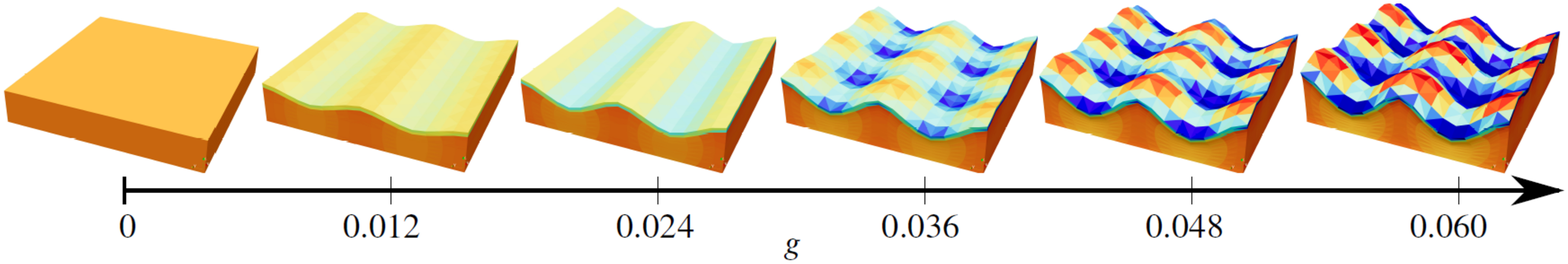}
\label{Fig8a_mu750}
\end{figure}
\begin{figure}[!h]
\raggedright
\rotatebox{90} {(e) $\mu_{fiber}=1000$}
\def\svgwidth{0.9\textwidth}
\centering
\includegraphics*[width=0.90\textwidth]{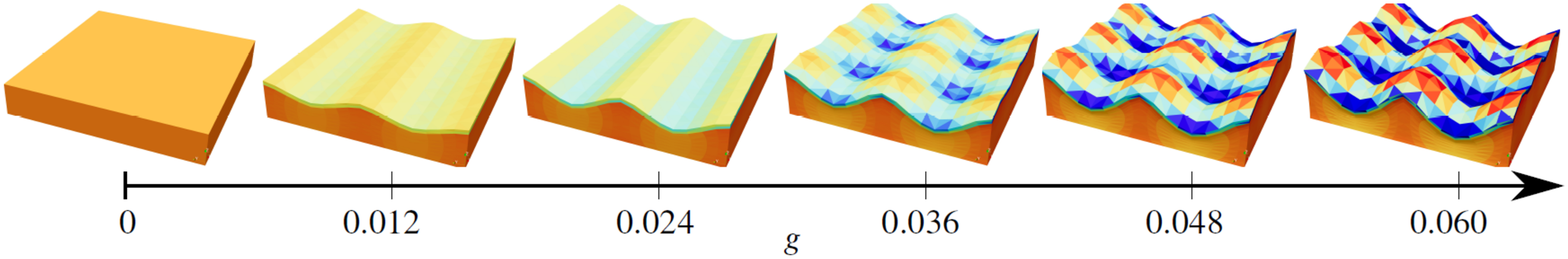}
\label{Fig8a_mu1000}
\end{figure}
\begin{figure}[!h]
\raggedright
\rotatebox{90} {(f) $\mu_{fiber}=2500$}
\def\svgwidth{0.9\textwidth}
\centering
\includegraphics*[width=0.90\textwidth]{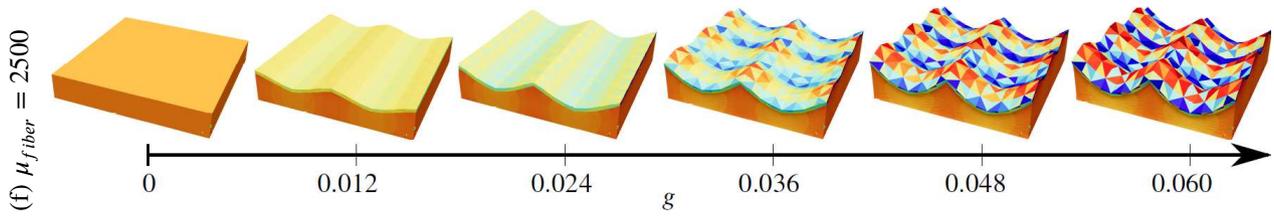}
\vspace{10pt}%
\caption{Pressure distribution of the buckling regime of a bilayer plate with monotonically increasing growth, $g$}
\label{Fig9}
\end{figure}
Figure \ref{fig10} presents the effect of fiber stiffness on the critical growth parameter that initiates buckling. While the vertical axis represents the fiber stiffness, the horizontal axis shows the planar growth parameter. Each deformed bilayer plate image corresponds to a state having a critical growth parameter that initiates either the first or the second buckling for each fiber stiffness. It is observed that the critical growth parameter decreases with the fiber stiffness for the first buckling. However, the secondary buckling modes are triggered within a small range of growth parameters close to each other. Since fibers are aligned to the y-axis, the fiber stiffness directly affects the critical growth in this direction. On the other hand, in the x-direction, which is perpendicular to fibers, the fibers do not affect material behavior. Nevertheless, they still are not fully decoupled from each other. The geometric form of the first buckling mode have also an effect on the secondary buckling at the post-buckling stage by shifting it to higher critical growth compared to isotropic (without fiber) case.
\begin{figure}[thb!]
\centering
\psfrag{mufiber}     [l][l] {\large{$\mu_{fiber}$}}
\psfrag{gg}     [l][l] {\large{$g$}}
\def\svgwidth{1.0\textwidth}
\includegraphics*[width=0.7\textwidth]{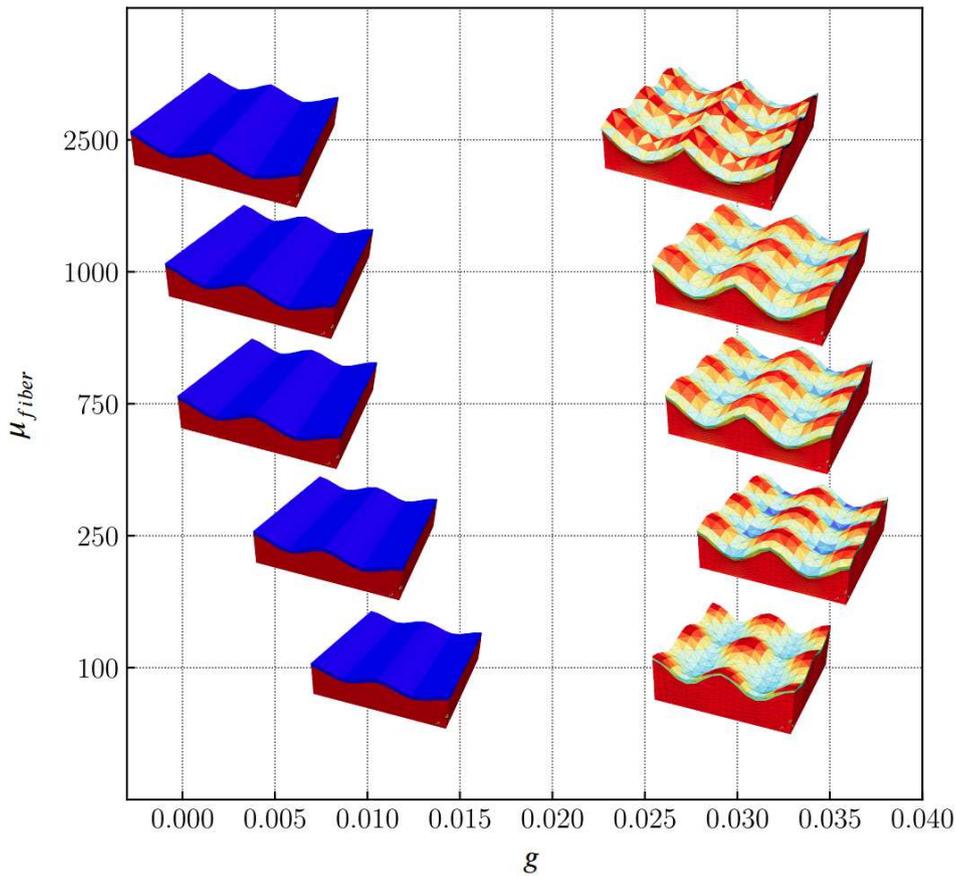}
\caption{Variation of the critical growth parameters for the first two bucking mode with respect to fiber stiffness.}
\label{fig10}
\end{figure}
Out-of-plane displacement and stored energy variations during the growth are also significant indicators to determine the critical growth parameter. The right column in Figures \ref{fig11} and \ref{fig12} shows the isotropic, volumetric, and anisotropic energy variation of the film and substrate separately. The left column represents vertical displacement change during growth at points A, B, and C. These specific points are located at the maximum and minimum displacement points on concave and convex regions of deformed shape. For each case in Figures \ref{fig11} and \ref{fig12}, it is clearly seen that vertical displacements of three points are the same until it reaches the bifurcation point where buckling initiates. In the isotropic bilayer plate, the vertical displacement patterns of the three points seem different from anisotropic ones. The reason is due to the isotropic nature of the bilayer plate, where the first buckling is simultaneously triggered both in x and y-directions at the same time. Later, the buckling shape starts to evolve to the labyrinth form at the second bifurcation point. Fiber-reinforced stiff film models show two different bifurcation points until 0.04 $g$ level. Thesebifurcations correspond to the first buckling initiation along the fiber direction, and the second buckling initiates perpendicular to the fiber direction. In alignment with Figure \ref{fig10}, it is observed that by increasing the fiber stiffness, the critical growth parameter decreases. For example, the first buckling initiates at $g=0.0116$ for $\mu_{fiber}=100$, while it initiates at $g=0.0026$ for $\mu_{fiber}=2500$. For the intermediate values of fiber stiffness, the critical growth parameter is revealed between this range. The second buckling initiation is observed in a small range of $g=0.028-0.032$ levels. The second buckling is affected by the geometrical shape of the first buckling mode. Different parts of the stored energies in the film and the substrate are shown in the right column in Figures \ref{fig11} and \ref{fig12}. Due to the sudden energy release, abrupt changes and kinks in energy plots are observed. Tracking energy is a better indicator than the displacement variations to identify the bifurcation as energy is scalar quantity and location independent. While the slope of the isotropic energy of the film shows a decrease at the first buckling for lower fiber stiffnesses, it increases for higher fiber stiffnesses. It is concluded that the first buckling is the combination of energy release of stiff film and stiff fibers. When the fiber stiffness increases, the energy contribution of fibers to the first instability increases as well. Furthermore, the isotropic energy of the film and substrate layer has a key role in the initiation of the secondary buckling.
\begin{figure}[thb!] 
\centering
\psfrag{NOF}     [l][l] {$\underline{no \  fiber}$}
\psfrag{mu100}     [l][l]  {\normalsize{$\underline{\mu_{fiber}=100}$}}
\psfrag{mu250}     [l][l] {\normalsize{$\underline{\mu_{fiber}=250}$}}
\psfrag{mu750}     [l][l] {\normalsize{$\underline{\mu_{fiber}=750}$}}
\psfrag{mu1000}     [l][l] {\normalsize{$\underline{\mu_{fiber}=1000}$}}
\psfrag{mu2500}     [l][l] {\normalsize{$\underline{\mu_{fiber}=2500}$}}
\def\svgwidth{1.0\textwidth}
\psfrag{yy1}     [l][l] {\normalsize{$Vertical \ displacement \ [-]$}}
\psfrag{yy2}     [l][l] {\normalsize{$Energy \ [-]$}}
\psfrag{gg}     [l][l] {\normalsize{$g$}}
\includegraphics*[width=0.85\textwidth]{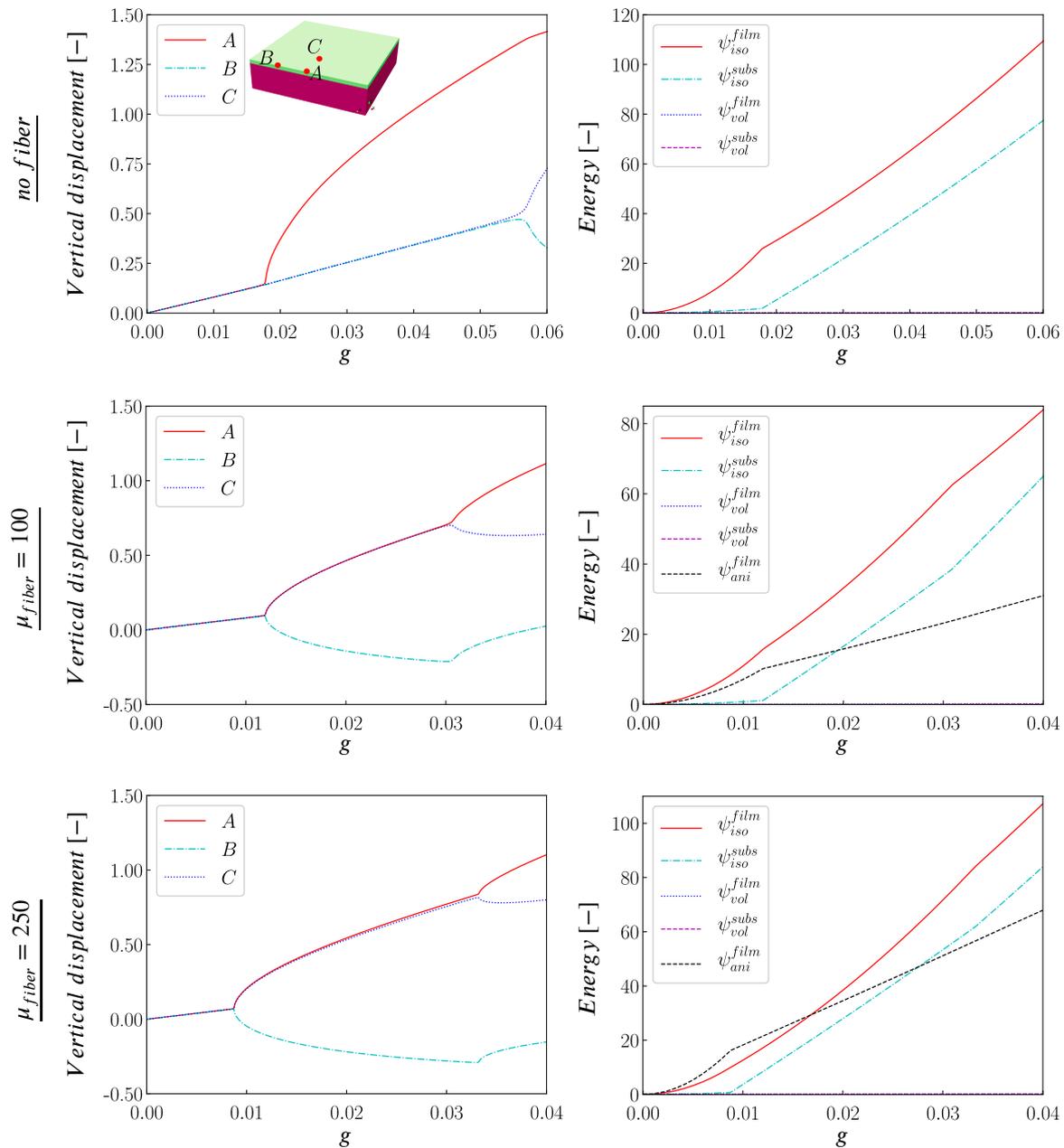}
\caption{Variation of vertical displacements of points A, B, and C; and variation of decoupled forms of energy for the film and substrate, those are subjected to incrementally growth for the following cases; without fiber reinforcement, $\mu_{fiber}$=100 and $\mu_{fiber}$=250, respectively.}
\label{fig11}
\end{figure}
\begin{figure}[thb!] 
\centering
\psfrag{mu100}     [l][l] {\normalsize{$\underline{\mu_{fiber}=100}$}}
\psfrag{mu250}     [l][l] {\normalsize{$\underline{\mu_{fiber}=250}$}}
\psfrag{mu750}     [l][l] {\normalsize{$\underline{\mu_{fiber}=750}$}}
\psfrag{mu1000}     [l][l] {\normalsize{$\underline{\mu_{fiber}=1000}$}}
\psfrag{mu2500}     [l][l] {\normalsize{$\underline{\mu_{fiber}=2500}$}}
\def\svgwidth{1.0\textwidth}
\psfrag{yy1}     [l][l] {\normalsize{$Vertical \ displacement \ [-]$}} 
\psfrag{yy2}     [l][l] {\normalsize{$Energy \ [-]$}}
\psfrag{gg}     [l][l] {\normalsize{$g$}}
\includegraphics*[width=0.85\textwidth]{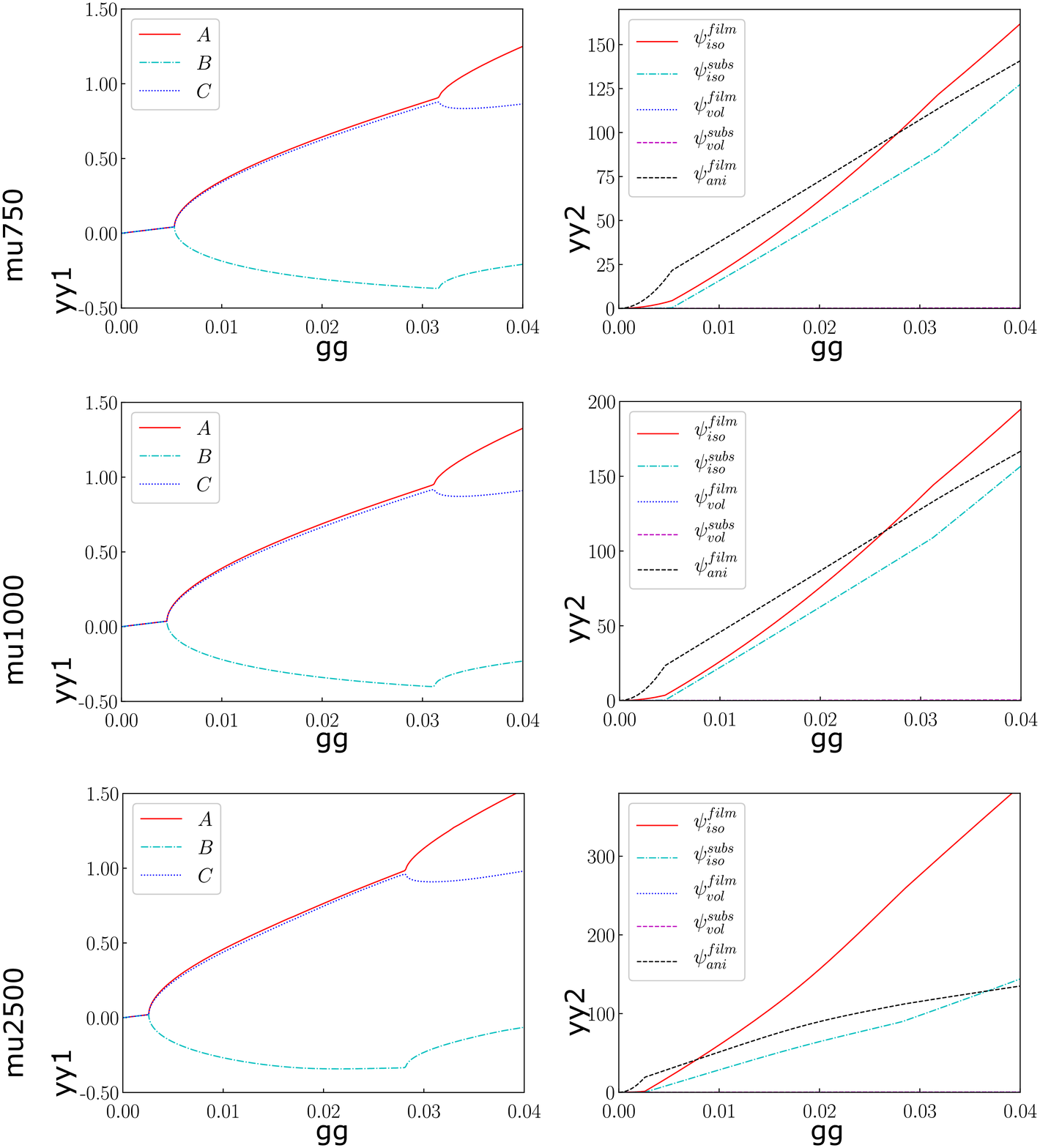}
\caption{Continuation of Figure \ref{fig11} - Variation of vertical displacements of points A, B, and C; and variation of decoupled forms of energy for the transversely isotropic film and isotropic substrate, those are subjected to incrementally growth for the following cases; $\mu_{fiber}$=750, $\mu_{fiber}$=1000, $\mu_{fiber}$=2500, respectively.}
\label{fig12}
\end{figure}
\section{Remarks and conclusion}\label{sec05}
In this work, planar growth-induced instabilities on bilayer 3D thick stiff film on the compliant substrate are examined in a various range of fiber stiffness with a five field Hu-Washizu type mixed variational formulation on Eulerian configuration for the quasi-incompressible and quasi-inextensible limits of transversely isotropic materials. We have studied a numerical example of mixed variational formulation by the implementation of \emph{T2P0F0} element, in order to understand the effect of fiber stiffness on a planar growth-induced primary and secondary instabilities on 3D representative bilayer plate. As a preliminary study for the post-buckling analysis, it was worked on a long but thin bilayer plate for the determination of exact wavelengths which plays key role in the periodic boundary condition in 3D plate having rectangular in-plane section. It is shown that, the wavelength decreases by increasing the fiber stiffness. According to results of second phase of the study, while higher fiber stiffness on the film layer causes the first instability in the direction of fibers with a lower growth parameter $g$, the effect of the fiber stiffness is in the minor level on the secondary buckling mode where it is observed perpendicular to fiber direction. Another outcome is the energy release mechanism at the initiation of buckling is mainly composed of isotropic and anisotropic contributions of stiff film layer. For a higher fiber stiffness, the effect of the anisotropic energy on the first buckling becomes more dominant than other type of energy components. However in the secondary instability, isotropic energy of the film layer plays key role for energy release mechanism. Numerical outcomes of this study will help to comprehend the role of the fiber stiffness  on the the buckling and post-buckling behavior of multi-layer structures. As a future work, the investigation of the growth-induced buckling on 3D biological tissues will be studied.

\section{Ethical statement}\label{sec0x}
None.


\end{document}